\documentclass[11pt,leqno]{article}
\usepackage{graphicx}
\usepackage{amsfonts, amsthm}
\usepackage{amsxtra}
\usepackage[latin1]{inputenc}
\usepackage{fontenc}
\usepackage[dvips]{epsfig}
\begin{document}
\newtheorem{theo}{Theorem}[section]
\newtheorem{defi}[theo]{Definition}
\newtheorem{lemm}[theo]{Lemma}
\newtheorem{prop}[theo]{Proposition}
\newtheorem{rem}[theo]{Remark}
\newtheorem{exam}[theo]{Example}
\newtheorem{cor}[theo]{Corollary}
\newcommand{\mat}[4]{
     \begin{pmatrix}
            #1 & #2 \\
            #3 & #4
       \end{pmatrix}
    }
\def\Z{\mathbb{Z}} 
\def\R{\mathcal{R}} 
\def\I{\mathcal{I}} 
\def\C{\mathbb{C}}   
\def\N{\mathbb{N}}   
\def\PP{\mathbb{P}}   
\def\Q{\mathbb{Q}}   
\def\L{\mathcal{L}}   
\def\ol{\overline}   

\newcommand{\gcD}{\mathrm {\ gcd}}  
\newcommand{\End}{\mathrm {End}}  
\newcommand{\Aut}{\mathrm {Aut}}  
\newcommand{\GL}{\mathrm {GL}}  
\newcommand{\SL}{\mathrm {SL}}  
\newcommand{\PSL}{\mathrm {PSL}}  
\newcommand{\Mat}{\mathrm {Mat}}  

\newcommand\ga[1]{\overline{\Gamma}_0(#1)}      
\newcommand\pro[1]{\mathbb{P}^1(\mathbb{Z}_{#1})}  
\newcommand\Zn[1]{\mathbb{Z}_{#1}}
\newcommand\equi[1]{\stackrel{#1}{\equiv}}
\newcommand\pai[2]{[#1:#2]}    
\newcommand\modulo[2]{[#1]_#2} 
\newcommand\sah[1]{\lceil#1\rceil}       
\def\sol{\phi} 
\begin{center}
{\LARGE\bf Transfer operators for $\Gamma_{0}(n)$ and the Hecke operators for 
the period functions of 
$\PSL (2,\mathbb{Z})$} 
\footnote{ MSC: 11F25, 11F67, 37C40

Keywords: transfer operators, congruence subgroups, Lewis equation, Hecke operators, period functions
} \\
\vspace{.25in} {\large {\sc J. Hilgert\footnote{Institut f\"ur Mathematik, 
TU-Clausthal, 
D-38678 Clausthal-Zellerfeld. E-mail: hilgert@math.tu-clausthal.de}, 
D. Mayer\footnote{IHES, F-91440 Bures sur Yvette. On leave of absence from Institut f\"ur Theoretische Physik, TU-Clausthal, D-38678 
Clausthal-Zellerfeld. E-mail: mayer@ihes.fr or  dieter.mayer@tu-clausthal.de}, 
H. Movasati\footnote{Institut f\"ur Theoretische Physik, TU-Clausthal, D-38678 
Clausthal-Zellerfeld. E-mail: hossein.movasati@tu-clausthal.de}}}
\end{center}


\begin{abstract}
In this article we report on a surprising relation between the transfer 
operators for the congruence 
subgroups 
$\Gamma_{0}(n)$ and the Hecke operators on the space of period functions for 
the modular group 
$\PSL (2,\mathbb{Z})$.
 For this we study special eigenfunctions of the transfer operators with 
eigenvalues $\mp 1$, 
which are also solutions of the Lewis equations for the groups $\Gamma_{0}(n)$
 and which are determined by  eigenfunctions of the transfer operator for the 
modular group 
$\PSL (2,\mathbb{Z})$. In the language of the Atkin-Lehner theory of old and 
new forms one should 
hence call them old eigenfunctions or old solutions of Lewis equation. 
It turns out that the sum of the components of these old solutions for the 
group $\Gamma_{0}(n)$ 
determine for any $n$  a solution of the Lewis equation for the modular group 
and hence also an 
eigenfunction of the transfer operator for this group.

Our construction gives in this way linear operators in the space of period 
functions for the group 
$\PSL (2,\mathbb{Z})$. Indeed these operators  are just the Hecke operators 
for the period functions 
of the modular group derived previously by Zagier and M\"uhlenbruch  using the 
Eichler-Manin-Shimura 
correspondence between period polynomials and modular forms for the modular 
group.
\end{abstract}


\section{Introduction}

This paper has three main ingredients. The first is the transfer operator from 
statistical mechanics 
which plays an important role in the ergodic theory of
dynamical systems and especially in the theory of dynamical zeta functions 
(see \cite{may1}, \cite{ rue}). 
Here we are interested in the transfer operators for the geodesic flow on the 
surfaces $\Gamma/\mathbb{H}$ 
for $\Gamma$ any of the congruence subgroups $\Gamma_{0}(n)$. These operators 
have been introduced in 
\cite{cha}, \cite{cha1} in the study of Selbergs zeta function for these 
groups. 

The second ingredient are certain  functions holomorphic in the cut plane, 
introduced by J.~B.~Lewis in \cite{lew} in his study of the Maass wave forms 
for $\PSL (2,\mathbb{Z})$. 
They were later named period functions by Zagier (see \cite {zag1}) because of 
their
close relation to the classical period polynomials in the Eichler, Manin , 
Shimura theory of periods for 
cusp forms.  
Period functions for the modular group are solutions of the so called Lewis 
equation
\begin{equation}
\label{2feb03}
\sol(z)=\sol(z+1)+\lambda z^{-2s}\phi \left(1+\frac{1}{z}\right)
\end{equation}
with $\lambda=\pm1$, which fulfill certain  growth conditions at infinity 
depending on the weight $s$. 
When this weight 
 satisfies 
 $\Re(s)=\frac{1}{2}$,  these solutions are in 1-1 correspondence with the 
Maass cusp 
forms (see \cite{leza2}). 
There is a simple relation between the transfer operator for $\PSL 
(2,\mathbb{Z})$ and the period functions: they are just the eigenfunctions of 
this operator  with eigenvalue $\pm 1$ (see \cite{cha1}). 
When $s$ is  a negative integer $s=-n$   
the space of polynomial solutions of the Lewis equation is in 1-1 
correspondence  
with the space of period polynomials of cusp forms for $\PSL (2,\mathbb{Z})$ 
(see \cite{zag1}). The
Eichler-Shimura-Manin theory of periods however tells us that this space of
period polynomials modulo a certain one dimensional space is isomorphic 
to the direct sum of two copies of the space of  cusp forms of weight $2n+2$ 
in the half plane. 

The space of cusp forms is extensively studied in number theory and in 
particular we have the  Hecke algebra acting on it.

A Theorem by Choie  and Zagier (see \cite{chza} \S 3 Theorem 2) gives a
criterion to find an explicit realization of the corresponding Hecke operators 
when acting on the space of 
period polynomials or more generally period functions. Generalizing the 
description of Hecke operators 
for Maass wave forms by 
Manin in  \cite{man1}, Choie and Zagier found (see \cite{chza}, Theorem 3) 
an explicit form for these 
Hecke operators in the space of
period polynomials.
 Their matrices, however, from which the Hecke operators are constructed via 
the well known slash action of 
the group $\Mat_{n}(2,\mathbb{Z})$ on smooth functions, have negative entries 
and hence their action 
is defined only for entire weights.
 
The third important ingredient in our paper is a new realization of the Hecke
operators on period functions due to T. M\"uhlenbruch (to appear in his 
thesis, see
 \cite{mul}). He uses the  matrices in the set
$$
S_n:=\left\{\mat abcd \mid a>c \geq 0, \ d>b\geq 0,\  ad-bc=n\right\}
$$
such that the Hecke operators $T_n$ have the form
\begin{equation}
\label{madar}
T_n:=\sum_{A\in S_n} A
\end{equation}
acting on the period functions via the slash operator, and
where the sum is taken in the free abelian group generated by $S_n$.
All the matrices from the set $S_{n}$ have positive entries and so 
one can define the Hecke operators on period functions for any weight 
$s\in\C$.

The relation between the transfer operator and the period functions was  
considered originally only 
for the modular group $\PSL (2,\mathbb{Z})$. In this case the Lewis equation 
is a scalar equation 
for scalar functions and its solutions can be  related explicitly to the 
period functions of modular 
and Maass wave forms.
In order to extend this theory to more general Fuchsian groups
Chang and Mayer began in a series of papers (see
\cite{cha1} and its references) to investigate the transfer operator approach 
to congruence subgroups 
like $\Gamma_{0}(n)$, $ \Gamma^{0}(n)$ or $ \Gamma(n) $.
 This lead them to transfer operators acting in Banach spaces of vector valued 
holomorphic functions. 
The eigenfunctions of these operators then fulfill general Lewis  equations in 
vector spaces whose 
dimension is just the index in $\PSL (2,\mathbb{Z})$ of the corresponding 
subgroup.

In the present paper we discuss special solutions of these Lewis equations for 
the groups $\Gamma_{0}(n)$ 
which are in fact determined by the solutions of the Lewis equation for the 
modular group 
$\PSL (2,\mathbb{Z})$.
 Hence our construction is somehow reminiscent of the theory of Atkin and Lehner 
of old and new forms 
(see \cite{atk}). The exact connection will be discussed in a forthcoming 
paper.
 
To state our main results and to sketch the
content of each section we have to fix the notations used throughout the text. 
For each integer $n$ let $\Mat_n(2,\mathbb{Z})$ 
(resp. $\Mat_*(2,\mathbb{Z})$) be the set of $2\times 2$-matrices 
  with integer entries and determinant $n$ (resp. nonzero determinant) and
$\R_n:=\Z[{\Mat_n(2,\Z)}]$ (resp. $\R:=\Z[{\Mat_*(2,\Z)}]$) the set
of finite linear combinations (with coefficients in $\Z$) of
the elements of $\Mat_n(2,\Z)$ (resp. $\Mat_*(2,\Z)$). Note that 
$\R=\cup_{n\in \Z} \R_n$ and $\R_n\cdot\R_m\subset \R_{nm}$. 
By definition we have  
$$\GL(2,\Z)=\Mat_1(2,\Z)\cup \Mat_{-1}(2,\Z).$$
The following four elements of $\GL(2,\Z)$ will play a prominent role in this 
paper:
$$
I:=\mat1001, \quad M:=\mat0110 , \quad  T:=\mat1101 , \quad Q:=\mat0{-1}10.
$$
It turns out that instead of the groups $\Gamma_{0}(n)$ it is more convenient 
to  use their extensions 
$\ga n$  in $\GL(2,\mathbb{Z})$
$$
\ga n:=\left\{\mat abcd \in \GL(2,\Z)\mid c\equiv 0 \mod n \right\}
=\Gamma_{0}(n)\cup \Gamma_{0}(n)\mat100{-1}. 
$$
In \S ~\ref{mul} we recall the definition of the transfer operators  for 
$\Gamma_{0}(n)$ and $\ga n$ 
as used by Chang and Mayer in \cite{cha}, \cite{cha1},  respectively by Manin 
and Marcolli in \cite{man}, 
discuss briefly their relation and derive the Lewis equation for the 
eigenfunctions of the operator of 
Manin and Marcolli. This operator is defined on a space of holomorphic 
functions with values in the 
representation space of $\GL(2,\Z)$
induced from the trivial representation of $\ga n$. In order to describe and solve the 
corresponding Lewis 
equation it turns out that the appropriate indexing of the  components of these 
functions by the set  
 $$
I_n:= \ga n\backslash \GL(2,\Z)
$$ helps a lot. The group 
$\GL(2,\Z)$ acts on this coset space on the right in a canonical way. 
The detailed structure of $I_n$ will be studied in \S~\ref{ind} (in particular 
see 
Propositions~\ref{heute} and \ref{hilg2}). 
The different components  of the Lewis equation can then be written for $i \in 
I_{n}$ as follows

\begin{equation}
\label{DFG}
\sol_i(z)-\sol_{(iT^{-1})}(z+1)-\lambda 
z^{-2s}\sol_{(iT^{-1}M)}\left(1+\frac{1}{z}\right)=0.
\end{equation}

These equations have to be solved simultaneously  with functions $\sol_i\ $  
holomorphic in  the 
cut plane $\C\setminus(-\infty,0]$ for all $ i\in I_n$.

Let $\I^\lambda:=(I-T-\lambda TM)\R$ be the right ideal
generated by $(I-T-\lambda TM)$ in $\R$
. Consider then the following system of 
equations 
in the right $\R$-module $\I^\lambda\backslash\R$
\begin{equation}
\label{DFG1}
\psi_i-\psi_{(iT^{-1})}T-\lambda\psi_{(iT^{-1}M)}TM=0\mod 
\I^\lambda,\quad\forall  i\in I_n,
\end{equation}
 which obviously is closely related to (\ref{DFG}).
Here the $\psi_i$'s are unknown elements in $\R$. Note the two
different matrix actions in these equations: on the one hand matrices acting 
from the right on the 
index 
$i$ of $\psi_i$ and on the  other hand matrices acting from the right on 
elements 
$\psi_i$ via the ring multiplication of $\R$. Moreover, in  (\ref{DFG}) we
have the familiar slash operation formally defined for $s \in \mathbb{C}$ and 
$R\in \R$  by 

\begin{equation}\label{slashdef}
\sol\mid_s R (z)= | \det R|^{s}(cz+d)^{-2s}\sol\left(\frac{az+b}{cz+d}\right).
\end{equation}

Now suppose $\psi_i, i\in I_n$ solves (\ref{DFG1}) and 
$\sol$ is a solution of the Lewis equation (\ref{2feb03}) for $\PSL 
(2,\mathbb{Z})$ . 
For $s$ an integer the left hand side of (\ref{DFG1}) can
act on $\sol$ via the usual  slash-operator 
and one
obtains a solution $(\phi_i)_{i\in I_n}$ of (\ref{DFG}) by setting
$$\phi_i:=\phi\mid_s \psi_i$$
since $\sol\mid_s \I^\lambda=0$. 

It is well known  that $\frac{1}{z}$ is up to a constant factor the only 
solution of the scalar Lewis 
equation (\ref{2feb03})
for $\lambda=1$ and $s=1$ (see \cite{Ma91}). It follows from a result by Y. 
Manin and M. Marcolli 
(see Proposition~\ref{30jan03} and 
Remark \ref{scalLewisinRIn})
  that for these parameter values $(\phi_i)_{i\in\I_n}$
with $\phi_i(z)=\frac{1}{z}$ for all $i\in I_n$ is, up to a trivial scalar 
factor, also the unique 
solution of (\ref{DFG}).
Hence, if $(\psi_i)_{i\in I_n}$ solves (\ref{DFG1}), then there exists a 
constant $\kappa$ such that
$$
\frac{1}{z}\mid_1 {\psi_i}=\kappa\frac{1}{z} \quad\forall  i\in I_n 
$$ must hold.
Suppose furthermore that $\psi_i= \sum_{A\in P_i} A$, 
where $P_i$ is some finite subset of  $\Mat_{n}(2,\Z)$. Then the above 
equality reads
$\sum_{A\in P_i}\frac{1}{(az+b)(cz+d)}=\kappa\frac{1}{z}$, where $A=\mat 
abcd$.
The right hand side of this expression  obviously has a pole and a zero only
at $0$ and $\infty$. Hence other poles and zeroes of the left hand side 
must cancel. 
This means, however, that the matrices $A \in P_{i}$ have to be chosen in a 
very specific way.
Explicit calculations for the groups $\Gamma_{0}(n)$ for small $n$ lead us to 
an operator 
$K:A\mapsto K(A)$ which attaches to every matrix $A \in P_{i}$ another matrix  
$KA$ whose action just 
cancels the poles and zeros generated by the action of $A$. In all cases 
considered only a finite 
number of matrices A were necessary to get the correct pole and zero 
structure.
 We later found that an operator similar to $K$ was indeed used already by 
Choie and Zagier \cite{chza} 
and also by M\"uhlenbruch (see  \cite{mul}, Lemma 9) in their completely 
different derivation of the Hecke 
operators. 

The explicit form of the map $K$  is given as
\begin{eqnarray*}
K:S_n\setminus Y_n&\to& S_n\setminus X_n\\
 \mat abcd&\mapsto& T^{\sah{\frac{d}{b}}}Q\mat abcd=
 \mat {-c+\sah{\frac{d}{b}}a}{-d+\sah{\frac{d}{b}}b}ab,
\end{eqnarray*}
where 
$$ 
X_n:=\left\{ \mat c{a}0{\frac{n}{c}},\ c\mid n,\ 0\leq  a<\frac{n}{c} \right\}
,\quad
Y_n:=\left\{ \mat c0{a}{\frac{n}{c}},\ c\mid n,\ 0\leq a<c \right\}
$$ 
and where for a real $r$ we have denoted by  $\sah r$ 
the integer  satisfying  $\sah r-1<r\leq\sah r $. 
With the usual notation of Gauss brackets we obtain  $\sah r=-[-r]$.
The  inverse of $K$ is given by
$$
K^{-1}:\mat abcd \mapsto MT^{\sah{\frac{a}{c}}}QM\mat abcd=
 \mat cd{-a+\sah{\frac{a}{c}}c}
{-b+\sah{\frac{a}{c}}d}
$$
(see Proposition \ref{opKiso} for all this).
Borrowing terminology from algebraic geometry one may call $K$ a rational 
automorphism of $S_{n}$.  

To each index $i\in I_n$ we will attach a matrix (see Definition 
\ref{ABdef})
$$
A_i=\mat cb0{\frac{n}{c}},
$$
where $ c\geq 1$, $c\mid n$ and 
$0\leq b<\frac{n}{c}$ satisfy $\gcD(c,\frac{n}{c},b)=1$ . The numbers $c$ and 
$b$ are then uniquely 
determined by the index $i$.

Starting now with a matrix $A_i\in X_n$ we apply $K$ repeatedly until we
get an element of $Y_n$ where the iteration stops. Since $K$ is 
injective, 
two such chains of  elements in $S_n$ are either equal or disjoint. 
For $i\in I_n$ we denote by $k_{i}$ the number such that $K^jA_i$ is  well-defined for 
$j\leq k_i$ and 
$K^{k_i}A_i\in Y_n$ (see 
Definition \ref{kidef}).
Obviously each element in $X_n\cap Y_n$ forms a one-element chain so that 
$k_i=0$ for $A_i\in X_n\cap Y_n$.
Our main result then is

\begin{theo}
\label{haupt} 
The matrices 
$$
\psi_i=\sum_{j=0}^{k_i}K^j(A_i),\quad  i \in I_n
$$
determine a solution of equations ~(\ref{DFG1}).
Acting by these matrices through the slash operator  on a solution $\sol$ of 
the Lewis equation (\ref{2feb03}) for the 
group $\PSL (2,\mathbb{Z})$ with weight $s$ 
gives a solution of equation ~(\ref{DFG}) for the group $\overline{\Gamma}_{0}(n)$ with the same weight $s$.  
\end{theo}

In the second part of the theorem we used the fact that 
the slash operator with weight $s$ an arbitrary
complex number is indeed well defined for the elements of $\R$ defining the $\psi_{i}$.
Details will be discussed in
\S~\ref{act}. In particular, Lemma~\ref{2bahman81} gives a condition which 
ensures the equation 
$$\sol\mid_s{(I-T-\lambda TM)R}=\big(\sol\mid_s{(I-T-\lambda TM)}\big)\mid_s 
R=0$$
for $R\in\R$ to make sense and thus enables us to construct solutions of 
(\ref{DFG})
from solutions of (\ref{DFG1}) in the way explained above.

Theorem \ref{haupt}  shows that any solution $\sol$ of the scalar Lewis 
equation (\ref{2feb03}) for $PSL(2,\mathbb{Z})$
determines a solution $(\sol\mid_s{\psi_i})_{i\in I_n}$ of the system 
(\ref{DFG}) 
of Lewis equations 
corresponding to $\ga n$.
Since the sum of the components of any solution of ~(\ref{DFG}) is again a 
solution of
the scalar Lewis equation, this fact together with Theorem ~\ref{haupt}  
allows us to define a 
linear operator $\tilde T_n$ mapping the space of solutions of the 
scalar
Lewis equation for $\PSL (2,\mathbb{Z})$ to itself. For these operators we 
find

\begin{prop}
\label{prop1} 
The operators $\tilde T_n$ and the Hecke operators $T_n$ defined in 
~(\ref{madar}) are related through 
$$
T_n=\sum_{d^2\mid n}\mat d00d\tilde T_{\frac{n}{d^2}}.
$$ 
In particular they coincide if and only if $n$ is a 
product of distinct primes. 

\end{prop} 
Thereby we identified the matrices $T_{n}$ and $\tilde T_{\frac{n}{d^2}}$ with the operators they define via the slash action. 
The operators $\tilde T_n$ have been constructed through special solutions of 
the Lewis equation for the 
congruence subgroups $\overline{\Gamma}_{0}(n)$. There arises immediately the question if 
this is also the case for 
the Hecke operators $T_n$.
 Indeed, it turns out that also the operators $\tilde T_{\frac{n}{d^{2}}}$ 
appearing in the above 
Proposition \ref{prop1} can be related to special solutions of the Lewis 
equation for the group 
$\overline{\Gamma}_{0}(n)$:
 one shows quite generally that any solution of the Lewis equation for a 
group $\overline{\Gamma}_{0}(m)$ 
determines a solution of the corresponding  equation for the group 
$\overline{\Gamma}_{0}(ml)$ for arbitrary 
$l \in \mathbb{N}$. 
Its components are just copies of the former's components  (see Proposition 
9.8.)
 Taking then as the solution for the group $\overline{\Gamma}_{0}(m)$ the solution of 
Theorem 1.1 we get in this way a 
solution for the group $\overline{\Gamma}_{0}(ml)$. The sum of its components gives  just 
$\mu$-times the operator 
$\tilde T_m$ where $\mu$ is the index of $\overline{\Gamma}_{0}(ml)$ in $\overline{\Gamma}_{0}(m)$.
 This shows that also the operators
$\tilde T_{\frac{n}{d^{2}}}$ can be constructed from special solutions of the 
Lewis equation for the 
group $\overline{\Gamma}_{0}(n)$ and hence from special eigenfunctions of the transfer 
operator for this group.

Our results depend in a crucial way on a modified one-sided continued fraction 
expansion for rational numbers
and closely related partitions of $\R$ described in \S~\ref{mod}. 

The technical results about the slash-operation are provided in \S~\ref{act} 
and the transfer operators for $\Gamma_{0}(n)$ and $\overline{\Gamma}_{0}(n)$ are introduced in \S~\ref{mul}.
The indexing coset space $\ga n\backslash \GL(2,\Z)$ is studied in detail in 
\S~\ref{ind}. In \S~\ref{opK} we derive and discuss the operator $K$
and in \S~\ref{lew} we describe various versions of the Lewis 
equations for the groups $\Gamma_{0}(n)$ and $\overline{\Gamma}_{0}(n)$ and relate them to each other. This 
allows us to construct special 
solutions of these equations. 
 Finally, in \S~\ref{Heckeop} we show how our results lead to a completely new 
approach to the Hecke 
operators on the space of period functions for $\PSL (2,\mathbb{Z})$ which 
basically  only uses the 
transfer operators for the congruence subgroups $\Gamma_{0}(n)$ respectively $\overline{\Gamma}_{0}(n)$. Work on the 
extension of this approach 
to the Hecke operators also for other groups like the congruence subgroups is 
going on at the moment.


\section{A modified continued fraction expansion}
\label{mod}
\def\part{P}

This section is basically inspired by the work of M\"uhlenbruch in \cite{mul} 
adapted appropriately to our needs. M\"uhlenbruch introduces in 
\cite{mul} a modified continued fraction expansion for 
positive rational numbers
and attaches to each $x\in\Q^+$ a suitable chain of elements of $\R$ 
called a partition of $x$. To explain his construction
we begin by collecting some facts which are  
standard in the theory of continued fractions (see \cite{hw}).
Consider the finite continued fraction expansion of $x$
$$x=[a_0,a_1,\ldots,a_N]:=
a_0+\frac{1}{a_1+\frac{1}{a_2+\ldots+\frac{1}{a_N}}}$$
and put $\frac{p_n}{q_n}:=[a_0,a_1,\ldots,a_n]$ for $0\leq n\leq N$.
Then  $\gcd(p_n,q_n)=1$, $q_n\geq 0$, and the recursion formulas
\begin{eqnarray}\label{cfeden1}
p_n &=&a_np_{n-1}+p_{n-2}\\  \label{cfeden2}
q_n&=&a_nq_{n-1}+q_{n-2}
\end{eqnarray}
hold. In particular, we  have
\begin{equation}\label{qorder}
q_0\le q_1<\ldots< q_N.
\end{equation}
Moreover, the following equations 
\begin{equation}
\label{puff2}
\frac{p_0}{q_0}<\frac{p_2}{q_2}<\cdots\le x\le 
\cdots<\frac{p_3}{q_3}<\frac{p_1}{q_1}
\end{equation}
and 
\begin{equation}\label{det1}
p_nq_{n-1}-q_n p_{n-1}=(-1)^{n-1},\  p_nq_{n-2}-q_n p_{n-2}=(-1)^{n}a_n
\end{equation} hold.
We are going to fill the above sequence (\ref{puff2}) with more rational 
numbers.
We do that for the left hand side of the sequence, 
the case we later use. Assume that $n$ is even. The sequence of numbers
$$[a_0,a_1,\ldots,a_{n-1},t]=\frac{tp_{n-1}+p_{n-2}}{tq_{n-1}+q_{n-2}},
\quad \mbox{for } t=0,\ldots, a_n$$
is then strictly increasing from 
 $\frac{p_{n-2}}{q_{n-2}}$ to 
$\frac{p_n}{q_n}$. We insert these numbers into the left hand side of 
(\ref{puff2}) and obtain the 
longer sequence :
\begin{equation}
\label{29j}
\ldots<\frac{p_{n-2}}{q_{n-2}}<\frac{p_{n-1}+p_{n-2}}{q_{n-1}+q_{n-2}}
<\ldots<\frac{(a_n-1)p_{n-1}+p_{n-2}}{(a_n-1)q_{n-1}+q_{n-2}}<
\frac{a_np_{n-1}+p_{n-2}}{a_nq_{n-1}+q_{n-2}}=
\frac{p_n}{q_n}<\ldots
\end{equation}
Here we have used the convention $\frac{1}{a+\frac{1}{0}}=0$. 
If we denote the rational numbers $x_{j}$ in this sequence by 
$x_{j}=\frac{p_{j}'}{q_{j}'}$ then  
$\gcd(p_{j}',q_{j}')=1$ and two consecutive numbers  
$\frac{p_{j}'}{q_{j}'}<\frac{p_{j+1}'}{q_{j+1}'}$ satisfy 
$q_{j}'<q_{j+1}'$ and
\begin{equation}\label{pstrichdet1}
p_{j+1}'q_j'-p_{j}'q_{j+1}'=1,
\end{equation}
where the last equality is a consequence of (\ref{det1}).
Recall (see \cite{hw}) that for $x\in\Q^+$
there is a unique sequence $a_0,\ldots,a_n\in \N$ such that $a_n>1$ and
$x=[a_0,\ldots,a_n]$: If $x=[b_0,\ldots,b_{m-1},1]$ for  
$b_0,\ldots,b_{m-1}\in \N$, 
 then obviously $x=[b_0,\ldots,b_{m-1}+1]$, 
and hence
$m=n+1$ and $a_0=b_0,\ldots, a_{n-1}=b_{n-1}, a_n=b_n+1$.
This will be used in the  following definitions.

\begin{defi}{\rm
Given $x\in\Q^+$, the modified continued fraction expansion of $x$ is the 
sequence 
$x_j,\ j=0,1,\ldots $ recursively defined by: 
\begin{itemize}
\item $x_0:=x=[a_0,\ldots,a_N]$ with $a_N>1$.
\item If $x_{j-1}=[b_0,\ldots,b_m]$ with $b_m>1$, then
$$
 x_j:=   \begin{cases}  
                       [ b_0,b_1,\ldots,b_{m-1} ] & \mbox{ if }   2 \not\vert\  
m \\
                       [ b_0,b_1,\ldots,b_{m}-1 ] & \mbox{ if }  2\mid m.     
        \end{cases}
$$ 
\end{itemize}
If $x_{j-1}=0$, then $x_j=-\infty$ and the sequence stops.}
\hfill$\Box$
\end{defi}

Note that the length of the modified continued fraction expansion of 
$x=[a_0,\ldots,a_N]$ 
with $a_N>1$ is
not greater than $\sum_{i=1, \mbox{\tiny  even}}^N a_i$.

\begin{prop}\label{mcfevsminpart}
Let $x\in \Q^+$ and $x_0,x_1,\ldots,x_{k-1},x_k$ with $x_0=x$ and 
$x_k=-\infty$ be its modified continued fraction expansion. If 
$x_j=\frac{p_j}{q_j}$
with $\gcd(p_j,q_j)=1$ and $q_j\ge 0$, then we have
$p_{j-1}q_j-p_j q_{j-1}=1$ for $j=1,\ldots,k$ and
$q_0>q_1>\ldots> q_{k-1}>q_k=0$.
\end{prop}

\begin{proof}
Suppose that $x_{j-1}=\frac{p_{j-1}}{q_{j-1}}=[b_0,\ldots,b_m]$ with $b_m>1$. 
If $m$ is odd we have $x_j=\frac{p_j}{q_j}=[b_0,\ldots,b_{m-1}]$ and the 
relation
$p_{j-1}q_j-p_j q_{j-1}=1$
follows from (\ref{det1}) applied to the continued fraction
$[b_0,\ldots,b_m]$, whereas 
$q_{j-1}>q_j$ is a consequence of the inequalities in (\ref{qorder}) for 
$[b_0,\ldots,b_m]$.

In the case where $m$ is even the same calculation leading to 
(\ref{pstrichdet1}) 
can be used to derive $p_{j-1}q_j-p_j q_{j-1}=1$ from the recursion relations 
for the 
continued fraction $[b_0,\ldots,b_m]$. Here $q_{j-1}>q_j$ follows also from
the recursion relations for the continued fraction $[b_0,\ldots,b_m]$.
\end{proof}

\begin{defi}\label{partdef}{\rm
 A sequence 
$x_0,x_1,\ldots,x_{k-1},x_k$ of rational numbers  is called an {\it admissible 
sequence of length} 
$k+1$ if the following property holds:  
if $x_j=\frac{p_j}{q_j}$, where  $\gcd(p_j,q_j)=1$ and $q_j\geq 0$, then 
\begin{equation}
\label{9jan03}
\det\mat{q_{j-1}}{-p_{j-1}}{q_{j}}{-p_{j}}=1 \quad \forall j=1,2,\ldots,k.
\end{equation}

Let $x$ be a positive rational number. A {\it partition} $\part$ of $x$ is 
an admissible sequence $x_0,x_1,\ldots,x_{k-1},x_k$ with $x_0=x$ and 
$x_k=-\infty$. The number $k+1$ is called the {\it length} of the partition.
We use the convention $-\infty=\frac{-1}{0}, 0=\frac{0}{1}$. 
A partition $\part$ of $x$ is called a {\it minimal partition} if
\begin{equation}
\label{29jan03}
q_0>q_1>\ldots> q_{k-1}>q_k=0
\end{equation}
}\hfill$\Box$
\end{defi}

\begin{rem}\label{minpartrem}
{\rm
{}From (\ref{9jan03}) it follows that $p_{j-1}q_j>p_jq_{j-1}$ which implies
 $x_{j-1}>x_j$ for all $j=1,2,\ldots,k$. Moreover, (\ref{9jan03}) shows that 
the equation
$p_{j-1}q_j\equiv 1 \mod {q_{j-1}}$   has a unique 
solution  $q_j$ with $\ 0\leq q_j<q_{j-1}$. Therefore each $x\in\Q^+$ has a 
unique minimal 
partition, which we denote by $\part_x$.
According to Proposition \ref{mcfevsminpart}
the modified continued fraction expansion of $x\in \Q^+$ satisfies
(\ref{9jan03}) and (\ref{29jan03}). Therefore it agrees with 
the minimal partition $P_x$.
We will show in Proposition \ref{ezdev} that there is indeed no partition 
whose length is less than 
the length of the minimal partition  which  justifies the name minimal 
partition.
}
\hfill$\Box$
\end{rem}

Throughout this paper
we will use the notations introduced in  Definition \ref{partdef}.

\begin{rem}\label{gleicheq}
{\rm  Let $x=x_0,x_1,\ldots,x_{k-1},x_k$ be a partition of
$x\in \Q^+$ and $x_j=\frac{p_j}{q_j}$ with $\gcd(p_j,q_j)=1$ and $q_j\ge 0$.
\begin{enumerate}
\item[(i)] The equation $p_{k-1}q_k-p_k q_{k-1}=1$ implies that 
$-p_k=1=q_{k-1}$.
If the partition is minimal Remark \ref{minpartrem} and the construction of 
the modified 
continued fraction expansion of $x$ shows that in addition we have 
$p_{k-1}=0$.
\item[(ii)]  If
$q_{j-1}=q_j$ for some $j\in \{1,\ldots,k-1\}$, then (\ref{9jan03})
shows that $q_{j-1}=q_j=1$, i.e., $x_j=p_j=p_{j-1}-1=x_{j-1}-1$.
\end{enumerate}
}\hfill$\Box$
\end{rem}

For a partition $\part$ of $x$ of length $k+1$ given by 
$x_0,x_1,\cdots,x_{k-1},x_k$ 
and any index $1\leq
l\leq k-1$, a simple calculation shows that for
$x_j=\frac{p_j}{q_j}$ with $\gcd(p_j,q_j)=1$ for $j=0,1,\ldots,k$ the 
sequence
\begin{equation}
\label{13jan03}
 x_0,\ldots,x_{l-1},\frac{p_{l-1}+p_l}{q_{l-1}+q_l},
x_l,\ldots,x_{k-1},x_k
\end{equation}
defines a new longer partition $\part(l)$ of $x$. We call it a Farey extension 
of partition $\part$. 
One can also 
introduce the inverse of this construction: if a partition $\part$
contains a triple of the type 
$\frac{p_{l-1}}{q_{l-1}},\frac{p_{l-1}+p_l}{q_{l-1}+q_l},\frac{p_l}{q_{l}}$, 
then one can
delete $\frac{p_{l-1}+p_l}{q_{l-1}+q_l}$ and obtains in this way  a shorter 
partition 
$\check{ \part}(l)$ of $x$ called a Farey reduction of $\part$.

\begin{prop}
\label{ezdev}
Every partition $\part$ of a rational number $x\in \Q^+$ can be obtained from 
the 
minimal
partition $\part_{x}$ of $x$ by a  finite number of Farey extensions 
$\part(l)$. The minimal 
partition $\part_{x}$ can be derived from any partition $\part$ by a finite 
number of Farey 
reductions $\check{\part}(l)$.
\end{prop}

\begin{proof}
Given a partition $x_0,x_1,\ldots,x_{k-1},x_k$ 
of $x$ with $x_j=\frac{p_j}{q_j}$ and $\gcd(p_j,q_j)=1$
it is enough to prove that if the sequence  $(q_j)_{j=0,\ldots,k}$ 
is not decreasing,  then there exists a number $l\in\{1,\ldots,k-1\}$ such 
that
$$\frac{p_{l}}{q_{l}}=\frac{p_{l+1}+p_{l-1}}{q_{l+1}+q_{l-1}}.$$
Since $q_{k}=0$ there exists for $(q_j)_{j=0,\ldots,k}$ not strictly 
decreasing an index $l\in \{1,\ldots,k-1\}$ such 
that 
$$q_{l}>q_{l+1}\quad\mbox{but}\quad q_l\ge q_{l-1}.$$

If $q_l> q_{l-1}$, then the
triple $x_{l-1},x_l,x_{l+1}$ must be of the form $\frac{p_{l-1}}{q_{l-1}},
\frac{e+mp_{l-1}}{f+mq_{l-1}},\frac{e+(m-1)p_{l-1}}{f+(m-1)q_{l-1}}$, where
$m\in\N$ and $\frac{e}{f}$ is the unique rational number such that 
$p_{l-1}f-q_{l-1}e=1$
and $0\le e<q_{l-1}$.

If $q_l= q_{l-1}$, then Remark \ref{gleicheq} shows that $q_{l-1}=q_l=1$ and
$x_{l-1},x_l,x_{l+1}$ is of the form $\frac{p_{l-1}}{1},
\frac{p_{l-1}-1}{1},\frac{p_{l+1}}{q_{l+1}}$. But then (\ref{9jan03})
shows that $p_{l+1}=(p_{l-1}-1)q_{l+1}-1$ so that
$$\frac{p_{l-1}+p_{l+1}}{q_{l-1}+q_{l+1}}=p_{l-1}-1,$$
which implies the claim also in this case.
\end{proof}

\begin{lemm}
\label{16jan03}
Let $P_x = (x_0,x_1,\ldots,x_k)$ be the minimal partition of $x\in \Q^+$. If
$x_j=\frac{p_j'}{q_j'}$ with $\gcd(p_j',q_j')=1$ and $q_j'\ge 0$, then we have
\begin{enumerate}
\item[{\rm(i)}]
$x<\frac{p_{j-1}'-p_j'}{q_{j-1}'-q_j'}$ for $j=1,2,\ldots,k$.
\item[{\rm(ii)}]
$\sah 
{\frac{xq_{j+1}'-p_{j+1}}{xq_j'-p_j'}}=p_{j-1}'q_{j+1}'-p_{j+1}'q_{j-1}'$
for  $j=1,2,\ldots,k-1$.
\end{enumerate}
\end{lemm}

\begin{proof} 
\begin{enumerate}
\item[(i)]
Let $x_0=[a_0,\ldots,a_N]$ be the continued fraction expansion of $x$ with 
$a_N>1$. If 
$p_n$ and $q_n$ are the corresponding denumerators and denominators
defined by (\ref{cfeden1}) and (\ref{cfeden2}), then
Remark \ref{minpartrem} shows that the sequence
$\ldots 
>\frac{p_{j-1}'}{q_{j-1}'}>\frac{p_j'}{q_j'}>\frac{p_{j+1}'}{q_{j+1}'}>\ldots$
is the same as (\ref{29j}) which can also be rewritten as
\begin{equation}
\label{29j'}
\ldots<\frac{p_{n-2}}{q_{n-2}}
=\frac{p_n-a_np_{n-1}}{q_n-a_nq_{n-1}}
<\frac{p_{n}-(a_n-1)p_{n-1}}{q_{n}-(a_n-1)q_{n-1}}
<\ldots<\frac{p_{n}-p_{n-1}}{q_{n}-q_{n-1}}<
\frac{p_n}{q_n}<\ldots,
\end{equation}
where $n$ is even. For two consecutive elements
$$\frac{p_{j}'}{q_{j}'}=\frac{(k-1)p_{n-1}+p_{n-2}}{(k-1)q_{n-1}+q_{n-2}},
\quad\mbox{and}\quad
\frac{p_{j-1}'}{q_{j-1}'}=\frac{kp_{n-1}+p_{n-2}}{kq_{n-1}+q_{n-2}}$$
in (\ref{29j'}) we have
$$\frac{p_{j-1}'-p_j'}{q_{j-1}'-q_j'}=
\frac{(kp_{n-1}+p_{n-2})-((k-1)p_{n-1}+p_{n-2})}{(kq_{n-1}+q_{n-2})-((k-1)q_{n-1}+q_{n-2})}=
\frac{p_{n-1}}{q_{n-1}}$$
and since $n$ is even (\ref{puff2}) shows that this is larger than $x$. 

\item[(ii)] 
There are two possible forms for three consecutive elements in the sequence
(\ref{29j'}). The first is
\begin{equation}\label{triple1}
\frac{p_{j-1}'}{q_{j-1}'}=\frac{p_n-(k-1)p_{n-1}}{q_n-(k-1)q_{n-1}},\quad
  \frac{p_{j}'}{q_{j}'}=\frac{p_n-kp_{n-1}}{q_n-kq_{n-1}},\quad
  \frac{p_{j+1}'}{q_{j+1}'}=\frac{p_n-(k+1)p_{n-1}}{q_n-(k+1)q_{n-1}},
\end{equation}
where $k=1,2,\ldots,a_n-1$. Then, using (\ref{det1}) and $n$ even, we obtain
\begin{eqnarray*}
&&\hskip -2em p_{j-1}'q_{j+1}'-p_{j+1}'q_{j-1}'=\\
&=&(p_n-(k-1)p_{n-1})(q_n-(k+1)q_{n-1})-(p_n-(k+1)p_{n-1})(q_n-(k-1)q_{n-1})\\
&=&2.
\end{eqnarray*}
On the other hand,
using $k\ge 1$ and $p_{n-1}-xq_{n-1},xq_n-p_n>0$  (again recall that $n$ is 
even), we calculate
\begin{eqnarray*}
\sah {\textstyle\frac{xq_{j+1}'-p_{j+1}}{xq_j'-p_j'}}
&=&\sah{\textstyle\frac{x(q_n-(k+1)q_{n-1})-(p_n-(k+1)p_{n-1})}{x(q_n-kq_{n-1})-(p_n-kp_{n-1})}}\\
&=&1+\sah{\textstyle\frac{p_{n-1}-xq_{n-1}}{k(p_{n-1}-xq_{n-1})+xq_n-p_n}} \\
&=&2.
\end{eqnarray*}
Thus (ii) if proved for triples of the form (\ref{triple1}).

The second type of triples appearing in (\ref{29j'}) is
\begin{equation} \label{triple2}
\frac{p_{j-1}'}{q_{j-1}'}=\frac{p_n+p_{n-1}}{q_n+q_{n-1}},\quad
  \frac{p_{j}'}{q_{j}'}=\frac{p_n}{q_n},\quad
  \frac{p_{j+1}'}{q_{j+1}'}=\frac{p_n-p_{n+1}}{q_n-q_{n+1}}
\end{equation}
with even $n$. This time we have
\begin{eqnarray*}
p_{j-1}'q_{j+1}'-p_{j+1}'q_{j-1}'
&=&(p_n+p_{n+1})(q_n-q_{n-1})-(p_n-p_{n-1})(q_n+q_{n+1})\\
&=& a_{n+1}+2
\end{eqnarray*}
and
$$
\sah {\textstyle\frac{xq_{j+1}'-p_{j+1}}{xq_j'-p_j'}}
=\sah{\textstyle\frac{x(q_n-q_{n-1})-(p_n-p_{n-1})}{xq_n-p_n}}
=1+\sah{\textstyle\frac{p_{n-1}-xq_{n-1}}{xq_n-p_n}}.
$$
But an easy calculation again using ~(\ref{det1}) shows that
 
$\sah{\frac{p_{n-1}-xq_{n-1}}{xq_n-p_n}}=a_{n+1}+1$ if and only if 
$\frac{p_n+p_{n+1}}{q_n+q_{n+1}}\leq x<\frac{p_{n+1}}{q_{n+1}}$, which, 

according to (\ref{29j}), is indeed the case. 

\end{enumerate}

\end{proof}

\begin{defi}\label{mdef}{\rm
Consider an admissible sequence $\part=(x_0,\ldots,x_k)$ of $x_0=x\in \Q^+$ 
with $x_j=\frac{p_j}{q_j}$ such that  $\gcd(p_j,q_j)=1$ and $q_j\ge 0$. To 
this partition we attach the 
following element $m(\part)$
of $\Z[\R_1]=\Z[\SL(2,\Z)]$
\begin{equation}
m(\part)=\mat{q_{0}}{-p_{0}}{q_1}{-p_1}+\ldots+\mat{q_{l-1}}{-p_{l-1}}{q_l}{-p_l}+
\mat{q_{l}}{-p_{l}}{q_{l+1}}{-p_{l+1}}+\ldots+\mat{q_{k-1}}{-p_{k-1}}{q_{k}}{-p_{k}}.
\end{equation}
}
\hfill$\Box$
\end{defi}

Given two admissible sequences $P_1=(x_0,x_1,\ldots,x_k)$ and
$P_2=(y_0,y_1,\ldots,y_l)$ with $x_k=y_0$ we can define the {\it join}
\begin{equation}\label{joindef}
P_1\vee P_2=(x_0,x_1,\ldots,x_k,y_1,\ldots,y_l)
\end{equation}
of $P_1$ and $P_2$, which is again admissible. Note that in this case we have  
\begin{equation}\label{mjoin}
m(P_1\vee P_2)=m(P_1)+m(P_2).
\end{equation}

$\GL(2,\Z)$ acts on rational numbers from the left in the usual way:
$$
\mat abcdx=\frac{ax+b}{cx+d}.
$$
For the next lemma we will need the corresponding right action:
$$
x\mat abcd =\mat abcd^{-1}x=\frac{dx-b}{-cx+a}
$$
\begin{lemm}
\label{puff}
Let $P=(x_0,x_1,\ldots,x_k)$ be an admissible sequence and $A=\mat abcd
\in \GL(2,\Z)$ with
\begin{equation}
\label{27j2003}
\frac{a}{c}\geq x_i,\ i=0,1,2,\ldots,k
\end{equation}
(which for $c=0$ simply means $a>0$).
Then 
$$
P\cdot A:=
  \begin{cases}  
               (x_0A,x_1A,\ldots,x_{k-1}A,x_kA)& \mbox{for }   \det A=1 \\
               (x_kA,x_{k-1}A,\ldots x_1A,x_0A)& \mbox{for }  \det A=-1
   \end{cases}
$$
defines an admissible sequence with the property
$$
m(P)A=
  \begin{cases}   
               m(P\cdot A)& \mbox{for }   \det A=1 \\
               Mm(P\cdot A)& \mbox{for }  \det A=-1,
  \end{cases}
$$
where $(m(P),A)\mapsto m(P)A$ is the multiplication in $\R$.
\end{lemm}

\begin{proof}
Condition (\ref{27j2003}) implies that for $x_j=\frac{p_j}{q_j}$ with 
$\gcd(p_j,q_j)=1$
and $q_j\geq 0$, the number $x_jA=\frac{dp_j-bq_j}{aq_j-cp_j}\in \Q$ and 
$\gcd(dp_j-bq_j,aq_j-cp_j)$ =1, 
since 
$(r,s)\begin{pmatrix}p\\ q\end{pmatrix}=1$ implies
$\big((r,s)A\big)\big(A^{-1}\begin{pmatrix}p\\ q\end{pmatrix}\big)=1$.
Moreover, for $\det A=1$ the matrix
\begin{equation}\label{mPA=m(PA)}
\mat{aq_{j-1}-cp_{j-1}}{-dp_{j-1}+bq_{j-1}}{aq_j-cp_j}{-dp_j+bq_j}
=
\mat{q_{j-1}}{-p_{j-1}}{q_j}{-p_j}
\mat abcd 
\end{equation}
has determinant $1$, which  implies that $P\cdot A$ is indeed admissible.
The equality $m(P)A=m(P\cdot A)$ is immediate from (\ref{mPA=m(PA)}). 
The case $\det A=-1$ can be  treated similarly.
\end{proof}

To simplify the notation we define for $r\in \Z$ and $m\in\N$  the number
$(r)_m\in\{0,\ldots,m-1\}$
as 
$$(r)_m\equiv r\mod m.$$ 
Moreover, having fixed $n\in\N$ once and for all, we attach to each 
$i\in\{1,\ldots, 
n-1\}$ relative prime
to $n$ the number  $\hat{\imath} \in \{1,\ldots,n-1\}$ with
$i\hat \imath \equiv 1\mod n$.

\begin{lemm}\label{XTslemma}
Suppose that $c,c',i\in \N_0$ satisfy
$ c,c'\geq 1,\quad  c,c'\mid n$ and  $\gcd(c,c')=1=\gcd(i,n)$.
\begin{enumerate}
\item[{\rm(i)}] The matrix
$$
X:=\mat {(c \hat\imath  -c')_{\frac{n}{c'}}+c'}{c'}{\frac{n}{c'}}0 .
\mat c{(c'i)_{\frac{n}{c}}}0{\frac{n}{c}}^{-1}
$$
is contained in $\GL(2,\Z)$.
\item[{\rm(ii)}]
Set
$$
x:=\left\{\frac{c'i-c}{\frac{n}{c}}\right\},\quad
y:=\left\{\frac{c'i}{\frac{n}{c}}\right\},\quad 
z:=\left\{\frac{c \hat\imath -c'}{\frac{n}{c'}}\right\}.
$$
and
$$s:=\frac{((c'i-c)_{\frac{n}{c}}+c-(c'i)_{\frac{n}{c}})}{n/c}.$$
Then
$P_{z}\cdot X$ is an admissible sequence 
beginning with $y$ and ending in 
$x-s$. Moreover, the join
$(\part_{z}\cdot X)\vee (\part_x\cdot T^{s})$ is well-defined and a partition 
of 
$y$. 

\item[{\rm(iii)}]
\begin{eqnarray*}
&&\hskip -4em
m(\part_{x}) \mat c{(ic'-c)_{\frac{n}{c}}}0{\frac{n}{c}}T
+
M\,m(\part_{z})
\mat {c'}{(\hat \imath c-c')_{\frac{n}{c'}}}0{\frac{n}{c'}}TM=\\
&=&m\big((\part_{z}\cdot X)\vee (\part_x\cdot T^{s})\big)
\mat c{(ic')_{\frac{n}{c}}}0{\frac{n}{c}}.
\end{eqnarray*}
\end{enumerate}
\end{lemm}

\begin{proof}
\begin{enumerate}
\item[(i)]
Note first that
\begin{eqnarray*}
X
&=&\frac{1}{n} 
\mat {(c\hat \imath -c')_{\frac{n}{c'}}+c'}{c'}{\frac{n}{c'}}0 .
\mat {\frac{n}{c}}{-(c'i)_{\frac{n}{c}}}0c\\
&=&\mat 
{\frac{((c\hat \imath -c')_{\frac{n}{c'}}+c')\frac{n}{c}}{n}}
{\frac{-((c\hat \imath -c')_{\frac{n}{c'}}+c')(c'i)_{\frac{n}{c}}+cc'}{n}}
{\frac{\frac{n}{c}\frac{n}{c'}}{n}}
{\frac{-\frac{n}{c'}(c'i)_{\frac{n}{c}}}{n}}\\
&=&\mat 
{\frac{((c\hat \imath -c')_{\frac{n}{c'}}+c')}{c}}
{\frac{-((c\hat \imath -c')_{\frac{n}{c'}}+c')(c'i)_{\frac{n}{c}}+cc'}{n}}
{\frac{n}{cc'}}
{-\frac{(c'i)_{\frac{n}{c}}}{c'}}.
\end{eqnarray*}
Since $c$ and $c'$ are relative prime and divide $n$ we have $cc'\mid n$ so 
that
$c\mid\frac{n}{c'}$ and $c'\mid \frac{n}{c}$. But then
$$c\mid (c'+(\hat \imath c-c')_{\frac{n}{c'}}\quad\mbox{iff}\quad c\mid 
(c'+\hat \imath 
c-c')$$
and the latter is evident. Similarly $c'\mid (ic')_{\frac{n}{c}}$ 
reduces to $c'\mid ic'$ which is clear
and 
$$n\mid\big(-(c'+(\hat \imath 
c-c')_{\frac{n}{c'}})(ic')_{\frac{n}{c}}+cc'\big)$$
reduces to
$n\mid (-\hat \imath c ic'+cc')$ which again is evident.
Thus the entries of $X$ are all integral. Since $\det X=-1$ is immediate
from the definition of $X$, this proves the claim. 
 
\item[(ii)] 
Suppose that $P_z$ is given by $(z=z_0,\ldots,z_m)$. According to
Remark \ref{minpartrem} we have $z=z_0>z_1>\ldots>z_m$
Note that for $X$  condition ~(\ref{27j2003}) is satisfied for $P_z$. 
In fact, the number $\frac{a}{c}$ in  (\ref{27j2003}) is 
$$
\frac{\left(\frac{((c \hat \imath 
-c')_{\frac{n}{c'}}+c')}{c}\right)}{\left(\frac{n}{cc'}\right)}
=z+\frac{{c'}^2}{n}>z\ge z_j
\quad\forall j=0,\ldots,m.
$$
Now Lemma \ref{puff} shows that $P_z\cdot X$ exists and since 
$\det X=-1$ the first element of 
$P_{z}\cdot X$ is given by
$$
z_mX=-\infty X=
\frac
{-\frac{(c'i)_{\frac{n}{c}}}{c'}}{-\frac{n}{cc'}}=
\frac{(c'i)_{\frac{n}{c}}}{\frac{n}{c}}=y, 
$$
whereas the last element of $P_{z}\cdot X$ is given by
\begin{eqnarray*}
z_0X
&=&\frac{(c \hat \imath -c')_\frac{n}{c'}}{\frac{n}{c'}}
\mat 
{\frac{((c \hat \imath -c')_{\frac{n}{c'}}+c')}{c}}
{\frac{-((c\hat \imath -c')_{\frac{n}{c'}}+c')(c'i)_{\frac{n}{c}}+cc'}{n}}
{\frac{n}{cc'}}
{-\frac{(c'i)_{\frac{n}{c}}}{c'}}\\
&=&\frac{(c'i)_{\frac{n}{c}}-c}{\frac{n}{c}}\\
&=&\frac{(c'i-c)_{\frac{n}{n}}}{\frac{n}{c}}-
\frac{((c'i-c)_{\frac{n}{c}}+c-(c'i)_{\frac{n}{c}})}{\frac{n}{c}}\\
&=&x-s.
\end{eqnarray*}
Note that for $T^s$ the number $\frac{a}{c}$ in  (\ref{27j2003}) is   
$\frac{1}{0}=\infty$.
Since $\det T^s=1$  Lemma \ref{puff}  shows that $P_x\cdot T^s$ exists, starts 
with
$x T^s=x-s$ and ends at $-\infty T^s=-\infty$. Thus the join $(P_z\cdot 
X)\vee(P_x\cdot T^s)$
exists and is a partition of $y$.

\item[(iii)]
An elementary calculation shows 
$$
\mat c{(c'i-c)_{\frac{n}{c}}}0{\frac{n}{c}}T=T^s\mat
c{(c'i)_{\frac{n}{c}}}0{\frac{n}{c}}
$$
and using the formula for $X$ derived in (i) we also find
$$
\mat {c'}{(c\hat \imath -c')_{\frac{n}{c'}}}0{\frac{n}{c'}}TM=
\mat {(c\hat \imath -c')_{\frac{n}{c'}}+c'}{c'}{\frac{n}{c'}}0=
X \mat c{(c'i)_{\frac{n}{c}}}0{\frac{n}{c}}.
$$
Now using (ii), (\ref{mjoin}),  and Lemma \ref{puff} we calculate
\begin{eqnarray*}
&&\hskip -4em m(\part_{x}) \mat c{(ic'-c)_{\frac{n}{c}}}0{\frac{n}{c}}T
+
M\,  m(\part_{z})
\mat {c'}{(\hat \imath c-c')_{\frac{n}{c'}}}0{\frac{n}{c'}}TM\\
&=&[m(\part_{x})T^{s}+ M\, m(\part_{z})X]\mat 
c{(ic')_{\frac{n}{c}}}0{\frac{n}{c}}\\
&=&[m(\part_{x}T^s)+ m(\part_z X)]\mat c{(ic')_{\frac{n}{c}}}0{\frac{n}{c}}\\
&=&m\big((\part_{z}\cdot X)\vee (\part_x\cdot T^{s})\big)
\mat c{(ic')_{\frac{n}{c}}}0{\frac{n}{c}}.
\end{eqnarray*}
\end{enumerate}
\end{proof}


\section{The slash operator for complex weight $s$}
\label{act}

\def\G{\mathcal G}
\def\F{\mathcal F}
\def\DS{\mathcal{DS}} 

Let $\F$ be the set of functions $\sol$  holomorphic in the domain 
$\C\setminus(-\infty,r]$ for some $r=r_\phi$ which we call a branching point 
of $\phi$. Note that this does not rule out 
that $\sol$ extends to the point $r$ as a holomorphic function.
In $\F$ we have the usual addition and multiplication of functions. If
$\sol_1,\sol_2\in \F$, then one can find
$r_{\sol_1\sol_2},r_{\sol_1+\sol_2}$  such that 
$r_{\sol_1\sol_2},r_{\sol_1+\sol_2}\leq \max\{r_{\sol_1},r_{\sol_2} \}$. 
We fix the branch of $\log z$ in $\C\setminus(-\infty,0]$
which coincides with the ordinary logarithm on $(0,\infty)$ and set
$z^s:=e^{s\log z}$ for $z\in \C\setminus(-\infty,0]$ and $s\in \C$. 
For each matrix $\mat abcd\in\G$ with
$$
\G:=\left\{\mat abcd \in \Mat_*(2,\Z)\mid 
(c> 0 \hbox { or }( c=0\  \&\  a,d > 0 ))
 \right\}
$$ one has
 $(cz+d)^s\in\F$. If $c=0$ and
$\phi\in\F$, then also $\phi(\frac{az+b}{d})\in \F$.
Consider the subset $\DS$ of $\F\times \G$ consisting of those pairs 
$\left(\phi,\mat abcd\right)$ such that
there exists a branching point $r_\phi$ for $\phi$ with
\begin{equation}
\label{24nov2002}
a-  c r_\sol  >0
\end{equation}

\begin{prop}\label{DS1}
Fix $s\in \C$. Then the formula
\begin{equation}
(\sol\mid_s R) (z)= | \det 
R|^{s}(cz+d)^{-2s}\sol\left(\frac{az+b}{cz+d}\right)
\end{equation}
defines a map
\begin{eqnarray*}
\DS&\to& \F\\
(\phi,R)&\mapsto& \phi\mid_s R.
\end{eqnarray*}
If $r_\phi$ is a branching point for $\phi$ satisfying (\ref{24nov2002})
for $R=\mat abcd\in \G$, then 
$$\max\left\{\frac{dr_\phi-b}{a-cr_\phi}, -\frac{d}{c}\right\}$$
is a branching point for $\phi\mid_s R$, where we interprete $-\frac{d}{c}$ as 
$-\infty$ if $c=0$.
\end{prop}

\begin{proof}
Suppose that $z>-\frac{d}{c}$. Then 
$$\frac{az+b}{cz+d}> r_\phi\quad\Leftrightarrow\quad
 (a-cr_\phi)z > dr_\phi -b$$
and if (\ref{24nov2002}) is satisfied the last inequality is equivalent to 
$$z> \frac{dr_\phi - b}{a-cr_\phi}.$$
Now the claim is immediate.
\end{proof}

\begin{rem}{\rm
The {\it slash-operation} from Proposition \ref{DS1} can be extended by 
linearity
to the subset $\DS_\Z$ of $\F\times \Z[\G]$ consisting of those pairs 
$\left(\phi, \sum_{j=1}^m  n_jR_m\right)$ for which all $(\phi,R_j)\in\DS$.
In fact, suppose that (\ref{24nov2002}) is satisfied for $(\phi,R_j)$ with 
branching points
$r_{\phi,j}$ for $\phi$, then (\ref{24nov2002}) is satisfied for all 
$(\phi,R_j)$ 
with branching points $\min_j r_{\phi,j}$. 
}\hfill$\Box$
\end{rem}

\begin{prop} \label{DS2}
Suppose that $R_1,R_2, R_1R_2\in \G$ and
$(\phi,R_1), (\phi\mid_s R_1,R_2), (\phi, R_1R_2)\in \DS$. 
Then for each $s\in \C$
we have 
$$(\phi\mid_s R_1)\mid_s R_2= \phi\mid_s (R_1R_2).$$
\end{prop} 

\begin{proof} We argue by analytic continuation.
Note first that for $R_j=\mat {a_j}{b_j}{c_j}{d_j}$ we have
$R_1R_2=\mat{a_1a_2+b_1c_2}{a_1b_2+b_1d_2}{c_1a_2+d_1c_2}{c_1b_2+d_1d_2}$
and since $R_1,R_2, R_1R_2\in \G$ the functions
$$\big((c_1a_2+d_1c_2)z+(c_1b_2+d_1d_2)\big)^{-2s}$$
 and 
$$(c_2z+d_2)^{-2s}\left(c_1\left(\textstyle\frac{a_2 
z+b_2}{c_2z+d_2}\right)+d_1\right)^{-2s}
 =(c_2z+d_2)^{-2s}\left(\textstyle\frac{(c_1a_2+d_1c_2) 
z+(c_1b_2+d_1d_2)}{c_2z+d_2}\right)^{-2s}$$
are holomorphic on $\C\setminus(-\infty,0]$ and agree on $(0,\infty)$, hence 
agree everywhere.
But then
\begin{eqnarray*}
&&\hskip -1em \big((\phi\mid_s R_1)\mid_s R_2\big)(z)=\\
&=& |\det R_2|^s (c_2 z+d_2)^{-2s} (\phi\mid_s R_1)\left(\textstyle\frac{a_2 
z+b_2}{c_2z+d_2}\right)\\
&=& |\det R_2|^s (c_2 z+d_2)^{-2s}  |\det R_1|^s 
    \left(c_1 \left(\textstyle\frac{a_2 
z+b_2}{c_2z+d_2}\right)+d_1\right)^{-2s} 
    \phi\left(\textstyle\frac{a_1 \left(\frac{a_2 z+b_2}{c_2z+d_2}\right)+b_1}
                   {c_1\left(\frac{a_2 z+b_2}{c_2z+d_2}\right)+d_1}\right)\\
&=& |\det R_1R_2|^s (c_2 z+d_2)^{-2s} \left(c_1\left(\textstyle\frac{a_2 
z+b_2}{c_2z+d_2}\right)+d_1\right)^{-2s} 
\phi\left(\textstyle\frac{(a_1a_2+b_1c_2)z+(a_1b_2+b_1d_2)}{(c_1a_2+d_1c_2)z+(c_1b_2+d_1d_2)}\right)\\
&=&\big(\phi\mid_s (R_1R_2)\big)(z).
\end{eqnarray*}

\end{proof}

\begin{rem}\label{G+F0}{\rm 
Set $$
\G^+:=\left\{\mat abcd \in \G\mid 
a> 0; b,d \ge 0 ))
 \right\}
$$
and
$$\F_0:=\{\phi\in \F\mid 0 \mbox{ is a branching point of } \phi\}.$$
Then $\G^+$ is a multiplicative subsemigroup of $\Mat_*(2,\Z)$ and we have
$\F_0\times \G^+\subseteq \DS.$
Moreover the slash-operation $\mid_s$ induces a semigroup action 
$\F_0\times \G^+\to \F_0$.
In fact, given $R=\mat abcd\in\G^+$  and $\phi\in \F_0$, Proposition \ref{DS1} 
shows that 
$0\ge \max\left\{-\frac{b}{a},-\frac{d}{c}\right\}$ is a branching point for  
$\phi\mid_s R$.
Then Proposition \ref{DS2} implies the identity
$(\phi\mid_s R_1)\mid_s \R_2=\phi\mid_s (R_1R_2)$
for all $R_1,R_2\in \G^+$. Of course  we can extend the  action to 
$\Z[\G]\subset \R$ by linearity. Note, finally, that $I, T, TM$ and $MTM$ are 
contained in $\G^+$, 
but $M$ is not. 
}\hfill$\Box$
\end{rem}
\def\T{{\cal T}}
Let 
$$
\T:=\{\mat abcd \in\G\mid a>0 \}
$$
\begin{lemm}
\label{2bahman81}
For all $\phi\in \F_0$ and $P\in \Z[\T]$ the following equality is  
well-defined:
$$
\sol\mid_s{(I-T-\lambda TM)P}=\big(\sol\mid_s(I-T-\lambda TM)\big)\mid_sP
$$
\end{lemm}

\begin{proof} According to Remark \ref{G+F0} we have $\phi\mid_s (I-T-\lambda 
TM)\in \F_0$.
Therefore, in view of Proposition \ref{DS2}  it suffices to check that
$$TA=\mat {a+c}{b+d}cd , TMA=\mat {a+c}{b+d}ab, A\in \G$$ 
and all three satisfy (\ref{24nov2002}) with $r_\phi=0$. But this just means
$$
\frac{a+c}{c},\frac{a+c}{a},\frac{a}{c}>0,
$$ 
which follows immediately from the  hypothesis. Here, of course,  we 
interprete 
$\frac{a}{c}$ as $\infty$ for $c=0$. 
\end{proof}

\begin{rem}{\rm
\label{gula}
Let $A=\mat abcd$. If  $A\in\T$, then 
$$
TMA=\mat {a+c}{b+d}{a}{b}, \quad MTMA=\mat
{a}{b}{a+c}{b+d}\in\T,
$$ 
and if  $A\in \Mat_*(2,\Z^+\cup\{0\})$, then 
$$
MATM=\mat {c+d}{c}{a+b}{a},\quad 
ATM=\mat {a+b}{a}{c+d}{c}\in \T.
$$
We will need these facts in the proof of our main theorem. 
}\hfill$\Box$
\end{rem}

\section{Transfer operators for $\Gamma_{0}(n)$ and 
$\overline{\Gamma}_{0}(n)$}
\label{mul}

Let $W$ be a $\mu$-dimensional complex vector space and $A,B\in \Aut_\C(W)$.
We assume the isomorphisms $A^n\in \Aut_\C(W)$ to be uniformly bounded in 
$n\in\N$
w.r.t. one and hence any norm on $\Aut_\C(W)$.
Consider the Banach space $\mathcal{B}(D)$ of 
holomorphic functions in the disc $D=\{z \in \mathbb{C}:\  |z-1| < 
\frac{3}{2}\}$ 
which are continuous on $\overline{D}$ with the sup norm. Then the operator 
$
\L_s: \mathcal{B}(D)\otimes W\rightarrow\mathcal{B}(D)\otimes W
$
with
\begin{equation}
\label{nemikhanad}
\L_sf(z)=\sum_{n=1}^\infty (z+n)^{-2s}A^{n-1}Bf\left(\frac{1}{z+n}\right)
\end{equation}
is a nuclear operator for $\Re(2s)>1$ in this Banach space and $\L_s$
extends to a meromorphic family of nuclear operators in 
the whole
$s$-plane with possible poles of order one at the points $s=\frac{1-k}{2}$ 
with $k\in \N_0$. 
The proof follows the same line of arguments as in \cite{cha}. In 
fact, using the $k$-th Taylor polynomial of $f$ at $0$ we have:
\begin{equation}
\label{22nov02}
\L_sf(z)=
\L_{s+\frac{k+1}{2}}\tilde f(z)
+\sum_{i=0}^k\zeta_{A,B}(i+2s,z+1)\frac{f^i(0)}{i!},
\end{equation}
where 
$$\tilde 
f(z):=z^{-k-1}\big(f\left(z\right)-\sum_{i=0}^k\frac{f^i(0)}{i!}z^{i}\big)$$
and 
$$\zeta_{A,B}(a,b)=\sum_{n=0}^\infty \frac{A^{n-1}B}{(b+n)^a}$$ 
is  a kind of Hurwitz zeta function.
 The first term on the right hand side in expression (\ref{22nov02}) is 
holomorphic in 
$\Re(s)>\frac{1-(k+1)}{2}$ and the second term 
has poles of order one at $\frac{1-i}{2},\ i=0,1,\ldots, k$
(the proof of this last statement is as for the usual Hurwitz zeta function,  
\cite{lan}, Chapter XIV). 
This proves our assertion.  

By a direct calculation we have
$$
\L_sf(z)-(A\L_s)f(z+1)=(z+1)^{-2s}Bf\left(\frac{1}{z+1}\right)
$$
Therefore any eigenvector $f$ of $\L_s$ with eigenvalue $\lambda$
satisfies the following three term functional equation:
$$
\lambda(f(z)-Af(z+1))=(z+1)^{-2s}Bf\left(\frac{1}{z+1}\right).
$$
It is convenient to make the change of variable $z\mapsto z-1$
and introduce the new function $\Phi(z)=f(z-1)$. For $\lambda\not=0$ the above
equation then takes the form:
\begin{equation}
\label{18nov2002}
\Phi(z)-A\Phi(z+1)=\lambda^{-1}z^{-2s}B\Phi\left(1+\frac{1}{z}\right)
\end{equation}
Since $f$ is defined in the disk $D$, $\Phi$ is defined
in the shifted disk $\{z:\  |z-2| \leq \frac{3}{2}\}$. As in \cite{cha1} one 
shows that any 
eigenfunction $f$ of the operator $\L_s$ can be extended holomorphically to 
the entire 
complex plane $\mathbb{C}$ cut along the line $(-\infty,-1]$. Hence the 
corresponding function 
$\Phi(z)$ is holomorphic in $\mathbb{C}\setminus(-\infty,0]$.  In what follows 
we are interested in 
solutions of ~(\ref{18nov2002}) in the domain $\C\setminus(-\infty, 0]$
for the eigenvalues $\lambda=\pm 1$. 
 In the scalar case $\mu=1$ with $ A,B=I$ equation
~(\ref{18nov2002})   was
introduced by J. Lewis in \cite{lew}. The derivation of his equation via the 
transfer
operator appeared independently in \cite{may}. There one can also find the 
conditions under which a 
holomorphic solution of equation (\ref{18nov2002}) determines an eigenfunction 
of the transfer operator 
with eigenvalue $\lambda$. An interesting property of the solutions of 
equation (\ref{18nov2002})
is described by the following proposition:
 
\begin{prop}
If $\lambda=\pm 1$ and $(BA^{-1})^2=I$, then any solution of equation 
~(\ref{18nov2002})
in $\C\setminus(-\infty, 0]$ satisfies
\begin{equation}
\label{23nov02}
\Phi(z)=\lambda z^{-2s}BA^{-1}\Phi\left(\frac{1}{z}\right).
\end{equation}
\end{prop}

\begin{proof}
The domain $\C\setminus(-\infty, 0]$ is invariant under $z\mapsto 
\frac{1}{z}$. 
We insert $\frac{1}{z}$ in  ~(\ref{18nov2002}),  multiply it by 
$\lambda z^{-2s}BA^{-1}$ and then subtract the result from ~(\ref{18nov2002}). 
Using
the hypotheses we get the equality in (\ref{23nov02}).
\end{proof}

Of special interest for the following is the case $s=1$:
For this let us suppose that  $A$ and $B$ are two invertible real matrices 
with
non-negative entries which satisfy 
\begin{equation}
\label{mehdidarsnemikhanad}
A{\mathbb I}={\mathbb I}=B{\mathbb I},
\end{equation}
where  ${\mathbb I}$ is a $\mu$-dimensional vector with all components equal 
to $1$. 
This is for instance the case for $A$ and $B$ permutation matrices. 
Then the  vector $\Phi'=\Phi'(z)$ with all entries equal to 
$\frac{1}{z}$ is obviously a solution of ~(\ref{18nov2002}) with 
$\lambda=1$ and $s=1$. 

Generalizing the analogous result for the scalar case $\mu =1$ in \cite{Ma91} 
one has

\begin{prop}
\label{30jan03}
$\Phi'$ is up to a constant factor the unique solution of ~(\ref{18nov2002}) 
for $\lambda=1$ and 
$s=1$ in the Banach space $\mathcal{B}(D)\otimes W$  .  There does not exist 
any other solution of 
equation (\ref{18nov2002}) in this space for the parameter values $s=1$ and 
$\lambda$ with 
$\mid \lambda \mid=1$.
\end{prop}

\begin{proof} The proof is a straightforward adaption from \cite{may1}, 
Appendix C, 
and \cite{man}.
\end{proof}

{\bf Induced representations:}
Let $G$ be a group and $H$ be a subgroup of finite index
$\mu=[G:H]$ of $G$. For  each representation
$\chi: H\rightarrow \End(V)$ we consider the induced representation
$\chi_G:G\rightarrow \End(V_G)$, where
$$V_G:=\left\{f\colon G\to V\mid  
f(hg)=\chi(h)f(g)\quad \forall g\in G,h\in H \right\} $$
and the action of $G$ is given by
$$\big(\chi_G(g)f\big)(x)=f(xg)\quad\forall x,g\in G.$$
If $V=\C$ and  the initial representation is trivial,
the induced representation $\chi_G$ is the right regular representation 
$\rho\colon 
G\to \GL(\C^{H\backslash G})$.
In fact, in this case $V_G$ is the space of complex valued left $H$-invariant 
functions
on $G$ or, what is the same, complex valued functions on $H\backslash G$, and 
the action is by right 
translation in the argument. This also shows that we can view $\rho$ as a 
homomorphism
$G\to \GL(\Z^{H\backslash G})$. Moreover, for each $g\in G$ the operators 
$\rho(g)^n\in \End_\C(\C^{H\backslash G})$  are uniformly bounded in $n\in\N$.

\begin{rem}{\rm  One can
identify $V_G$ with $V^\mu$ using a set $\{g_1,g_2,\ldots,g_\mu\}$ of 
representatives 
for $H\backslash G$, i.e.
$$H\backslash G=\{Hg_1,Hg_2,\ldots,Hg_\mu \}.$$ 
Then 
\begin{eqnarray*}
V_G&\to& V^\mu\\
f&\mapsto& \big(f(g_1),\ldots,f(g_\mu)\big)
\end{eqnarray*}
is a linear isomorphism which transports $\chi_G$ to the linear $G$-action on 
$V^\mu$
given by
$$g\cdot(v_1,\ldots,v_\mu)
=\big(\chi(g_1gg_{k_1}^{-1})v_{k_1},\ldots,\chi(g_\mu 
gg_{k_\mu}^{-1})v_{k_\mu}\big),$$
where $k_j\in\{1,\ldots,\mu\}$ is the unique index such that $Hg_jg =H 
g_{k_j}$.
To see this one simply calculates
$$\chi_G(g)f(g_j)=f(g_jg)=f(g_jgg_{k_j}^{-1}g_{k_j})=\chi(g_jgg_{k_j}^{-1})
f(g_{k_j}).$$
In the case of the right regular representation
the identification $V_G\cong \C^\mu$ yields a matrix realization 
$$
\rho(g)=(\delta(g_igg_j^{-1}))_{i,j=1,\ldots,\mu},
$$
where $\delta(g)=1$ if $g\in H$ and $\delta(g)=0$ otherwise. 
Note for the following that the matrix $\rho(g)$ is a permutation matrix for 
all $g\in G$
\hfill$\Box$}
\end{rem}

In this article we are primarily intested in the subgroups $\Gamma_{0}(n) 
\subset \PSL (,2,\mathbb{Z})$, 
respectively their extensions  $\ga n \subset \GL(2,\Z)$. The representation 
$\chi$ is in both cases 
the  trivial 
representation of $\Gamma_{0}(n)$, respectively $\ga n$. 
The transfer operators for the groups $\Gamma_{0}(n)$ and $\ga n$ have been 
introduced by Chang and 
Mayer (see \cite{cha}, \cite{cha1}), respectively Manin and Marcolli (see 
\cite{man}).
 Taking in expression (\ref{nemikhanad}) for $A$ the matrix $\rho(QT^{\pm 
1}Q)$ and for $B$ the matrix 
$\rho(QT^{\pm 1})$ we get the transfer operators 
$\mathcal{L}_{s,\pm }$ for $\Gamma_{0}(n)$  whereas for $A=\rho(T^{-1})$ and 
$B=\rho(T^{-1}M)$ we have 
the transfer operator 
$\L_s $ for the group $\ga n$. 
An easy calculation shows that the operators 
$\mathcal{L}_{s,+}\mathcal{L}_{s,-}$ and $\L_s^{2}$ 
can be conjugated by the matrix $\rho( \mat100{-1})$.
On the other hand it was shown in \cite{cha} that the Selberg zeta function 
$Z_{\Gamma_{0}(n)}(s)$ for the group $\Gamma_{0}(n)$ can be expressed in terms 
of the Fredholm determinant 
of the operator $\mathcal{L}_{s,+} \mathcal{L}_{s,-}$ as
 $Z_{\Gamma_{0}(n)}(s)= \det(1-\mathcal{L}_{s,+}\mathcal{L}_{s,-})$ and hence 
also as 
$Z_{\Gamma_{0}(n)}(s)= 
\det(1-\mathcal{L}_{s}^{2})=det(1+\mathcal{L}_{s})\det(1-\mathcal{L}_{s})$. 
This shows that using the operator $\L_s$ the Selberg zeta function for the 
group $\Gamma_{0}(n)$ 
factorizes as in the case of the modular group and hence this transfer 
operator facilitates also the 
discussion of the period functions for $\Gamma_{0}(n)$. In the following we  
will 
therefore use this operator. 
The Lewis equation for $\Gamma_{0}(n)$ derived from the eigenfunction equation 
for $\mathcal{L}_{s}$ then has the form
\begin{equation}
\label{27JAN}
\Phi(z)-\rho(T^{-1})\Phi(z+1)-\lambda^{-1} z^{-2s}\rho(T^{-1}M)
\Phi\left(1+\frac{1}{z}\right)=0
\end{equation}

For the transfer operators considered above one finds $(BA^{-1})^2=I$ since 
$BA^{-1}=\rho(QTQT^{\pm1}Q)$,
 respectively
$BA^{-1}=\rho(T^{-1}MT)$,   and hence the two term equation ~(\ref{23nov02}) 
holds. 
Note that the  matrices in both examples are permutation matrices and so 
also the scalar equations in ~(\ref{23nov02}) involve only two terms.

\section{The indexing coset space}
\label{ind}

In this section we study the fine structure of $\ga n\backslash \GL(2,\Z)$
as a right $\GL(2,\Z)$-space. To do this we embed $\ga n\backslash \GL(2,\Z)$ 
into a natural 
$\GL(2,\Z)$-space with an action by a kind of linear fractional 
transformations.
We start with $\Z^2=\Z\times\Z$ on which $\GL(2,\Z)$ acts via
$$(x,y)\mat abcd=(ax+cy,bx+dy).$$
We define an equivalence relation $\sim_n$ on $\Z\times\Z$ via
$$(x,y)\sim_n(x',y')\quad:\Leftrightarrow\quad
(\exists k\in \Z) \gcD(k,n)=1, 
\begin{array}{c}
k x\equiv x' \mod n\\  
k y\equiv y' \mod n.
\end{array}
$$
Then the linearity of the action shows that it preserves $\sim_n$ so that the 
space
$[\Z:\Z]_n:=(\Z\times\Z)/\sim_n$ of equivalence classes inherits a
right $\GL(2,\Z)$-action. If, for fixed $n$, the equivalence class of $(x,y)$ 
is denoted by
$[x:y]$, then this action is given by
$$ [x:y]\mat abcd=[ax+cy : bx+dy]$$
which is of course very reminiscent of linear fractional transformations.
Note, however, that even for $n=p$ prime the space $[\Z\times\Z]_p$ is not the 
projective space $\PP^1(\Z_p)$ since we have not excluded the pairs of numbers 
both divisible by $p$.

\begin{rem}\label{Ininjection}
 {\rm The stabilizer of the point $[0:1]\in [\Z\times\Z]_n$ is $\ga n$ since
$(0,1)\sim_n (c,d)$ if and only if $c\equiv 0\mod n$ and $\gcD(d,n)=1$. Thus 
the orbit map
\begin{eqnarray*}
\GL(2,\Z)&\to& [\Z\times\Z]_n\\
g&\mapsto& [0:1] g
\end{eqnarray*}
factors to the equivariant injection
\begin{eqnarray*}\overline\pi\colon \ga n\backslash \GL(2,\Z)&\to& 
[\Z\times\Z]_n\\
\ga n \mat abcd&\mapsto& [c:d]
\end{eqnarray*}
}\hfill$\Box$

\end{rem}

Now we set $I_n:=\mathrm{Im}(\overline{\pi})\subseteq [\Z\times\Z]_n$ and note that $I_n$ 
is
$\GL(2,\Z)$-invariant.

\begin{prop}\label{gcd}
$I_n=\{[x:y]\mid \gcD(x,y,n)=1\}$.
\end{prop}

\begin{proof} ``$\supseteq$'': If $\gcD(x,y,n)=1$ set $m:=\gcD(x,y)$ and
$x':=\frac{x}{m}, y':=\frac{y}{m}$. Then $\gcD(m,n)=\gcD(x',y')=1$ and one can 
find 
$a,b\in \Z$ such that $ay'-bx'=1$. Therefore
$g:=\mat a{-b}{x'}{y'}\in \GL(2,\Z)$ and  
$$[0:1]g=[x':y']=[mx':my']=[x:y].$$
 
\noindent
``$\subseteq$'': If $[x:y]=[0:1]\mat abcd =[c:d]$, then there exist $k,r,s\in 
\Z$ such that
$\gcD(k,n)=1$ and 
\begin{eqnarray*}
k c -x&=&rn\\
k d -y&=&sn.
\end{eqnarray*}
If now $a=\gcD(x,y,n)$, then $a\mid \gcD(c,d)=1$  since $\mat abcd$ has 
determinant $ad-bc=1$.
\end{proof}

\begin{lemm} \label{Inrepr}
Given $m,n\in\Z$ and $u,v\in \Z$ such that $c=um+vn=\gcd(m,n)$ one can find 
$t\in \Z$ such that 
$$\gcd(u+\frac{n}{c}t,n)=1.$$
\end{lemm}

\begin{proof} Let $n=\prod_{j=1}^s p_j^{\alpha_j}$ be the decomposition into 
prime factors
and suppose that they are arranged in such a way that
$c=\prod_{j=1}^s p_j^{\beta_j}$ with
$\alpha_j=\beta_j$ for $j\le s_1$ and $\alpha_j>\beta_j$ for $j>s_1$.
Then $u\frac{m}{c}+v\frac{n}{c}=1$ implies that $u$ cannot contain a prime 
factor $p_j$ with $j>s_1$
so that
$\gcd(u,n)=\prod_{j=1}^{s_1} p_j^{\gamma_j}$
with $0\le \gamma_j\le \alpha_j$. We may assume w.l.o.g. that
$\gamma_j>0$ for $j\le s_2$ and $\gamma_j=0$ for $s_2<j\le s_1$, i.e.
$$\gcd(u,n)=\prod_{j=1}^{s_2} p_j^{\gamma_j}.$$
Now we pick $t=\prod_{j=s_2+1}^{s_1} p_j$ and comparing which $p_j$ divide 
respectively
$u,t,$ and $\frac{n}{c}$, we see that
no $p_j$ divides $u+\frac{n}{c}t$.
\end{proof}

\begin{prop}
\label{heute}
Each element of $I_n$ can be written as 
$\pai cd,\ c\geq 1,\ c\mid n$.
Here $c$ is determined uniquely, whereas $d$ is determined only up to an 
integer
multiple of $\frac{n}{c}$. It is possible to choose  $d=k d'$ with
$d'\geq 1, d'\mid n$, $1\le k< n$ and $\gcd(c,d')=1=\gcd(k,n)$.
\end{prop}

\begin{proof}
For $[x:y]\in I_n$ set $c=\gcd(x,n)$ and choose $u,v\in \Z$ such that 
$ux+vn=c$.
Using Lemma \ref{Inrepr} we can find  $t\in \Z$ such that
$\gcd(u+\frac{n}{c}t,n)=1$. Set $\nu:=u+\frac{n}{c}t$.
Then we have
$$\nu x = c + n \left(t\frac{x}{c}-v\right)\equiv c \mod n$$
so that with $d:=\nu y$ we obtain $[x:y]=[c:d]$. Now we set
$d':=\gcd(n,d)$ and use Lemma \ref{Inrepr}  in order to find a
$\nu'$ with $\gcd(\nu',n)=1$ such that
$\nu' d=d'\mod n$. Choosing $k\in \{1,\ldots,n-1\}$ such that $k\nu'\equiv 
1\mod 
n$
we have $d\equiv k d' \mod n$ and find
$[x:y]=[c:k d']$. This proves the existence part of the proposition
since $\gcd(c,d')=\gcd(x,y,n)=1$.

To prove  uniqueness suppose that $\pai cd=\pai {c'}{d'}$.
Then we have 
$c=lc'+rn,d=ld'+sn $ for some $r,s,l\in\Z$ with  $\gcd(l,n)=1$.
If now  $c,c'\geq 1,\ c,c'\mid n$ the first equality  implies that $c=c'$ and 
$l=-r\frac{n}{c}+1$.
Inserting this  $l$ into the second equality we obtain the uniqueness of $d$
up to $\frac{n}{c}\Z$.  
\end{proof}
Unfortunately the parametrization of the elements of $I_{n}$ by $[c:d]$ with 
$c \geq 1$ and $c \mid n$ 
is not unique as shown in Proposition \ref{heute}. To achive an unique 
parametrization we proceed as 
follows: 

\begin{defi}\label{parametrization}{\rm
For fixed $n \in \mathbb{N}$ and $c\in\{1,\ldots,n\}$ with $c \mid n$ choose 
$b\in \{0,\ldots,\frac{n}{c}-1 \}$.
 We call the pair $(c,b)$  {\it $n$-admissible} if there exists $k \in 
\{0,\ldots,c-1\}$ 
with $\gcd (c,b+k\frac{n}{c})=1$. For such a pair we set 
\begin{equation}\label{d(c,b)}
d(c,b):= min\{c+b+k\frac{n}{c}:k \in \{0,\ldots,c-1\}, 
\gcd(c,b+k\frac{n}{c})=1\}
\end{equation}
}\hfill$\Box$
\end{defi}

\begin{rem}\label{hilg}{\rm
\begin{enumerate}
\item[(a)] If $\gcd(c,b)=1 $, then $d(c,b)=c+b$.
\item[(b)] If $(c,b)$ is $n$-admissible, then $\gcd(c,b,\frac{n}{c})=1$.
\item[(c)] the pair $(c,b)$ is $n$-admissible iff $ \exists k \in \mathbb{Z}$ 
with 
$\gcd(c,b+k\frac{n}{c})=1$.
\end{enumerate}
}\hfill$\Box$
\end{rem}
We need the following lemma\footnote{We thank Ch.~Elsholtz for showing us how 
to prove this lemma}:
\begin{lemm}\label{hilg1}
Given the numbers $a$,$b$,$c \in \mathbb{Z}$ then $\gcd(a,b,c)=1$ iff there 
exists a $k\in \mathbb{Z}$ 
such that $\gcd(a,b+kc)=1$.
\end{lemm}

\begin{proof}
If $\gcd(a,b+kc)=1$ for some $k\in \mathbb{Z}$ then there exist $x,y \in 
\mathbb{Z}$ such that 
$ax+(b+kc)y=1$ and hence $\gcd(a,b,c)=1$.\\
Conversely, if $\gcd(a,b,c)=1$  define
$$ t_{ab}:=\gcd(a,b),\quad
t_{bc}:=\gcd(b,c),\quad
 t_{ac}:=\gcd(a,c),$$ respectively
 $$t_{a}:=\frac{a}{t_{ab}t_{ac}},\quad 
 t_{b}:=\frac{b}{t_{ab}t_{bc}},\quad 
t_{c}:=\frac{c}{t_{ac}t_{bc}}.$$
 Then 
$\gcd(t_{x},t_{y})=\gcd(t_{xy},t_{xz})=\gcd(t_{x},t_{yz})=1$ for all  $x\not=y\not=z \in \{a,b,c\}$.
Obviously 
$$a=t_{a}t_{ab}t_{ac},\ b=t_{b}t_{ab}t_{bc}, \ 
c=t_{c}t_{ac}t_{bc}.$$ 

To determine $k \in \mathbb{Z}$ with $\gcd(a,b+kc)=1$ we proceed as follows: 
In the case $\gcd(t_{ab},t_{a})=1$ one finds  for $a=t_{a}t_{ab}t_{ac}$ and 
$b+kc=t_{b}t_{ab}t_{bc}+kt_{c}t_{ac}t_{bc}$ with $k=t_{a}$  

\begin{eqnarray*}
1&=&\gcd(t_{ac},b)=\gcd(t_{ac}, b+kc),\\ 
1&=&\gcd(t_{ab},t_{a}t_{c}t_{ac}t_{bc})=\gcd(t_{ab}, b+kc),\\
1&=&\gcd(t_{a},b)=\gcd(t_{a}, b+kc).
\end{eqnarray*}
In the case $d=\gcd(t_{a},t_{ab})>1$ with $t_{a}=dt_{a}'$ and 
$t_{ab}=dt_{ab}'$ write $t_{a}'=sk$ with $s\mid d$ and $\gcd(k,d)=1$. Then $\gcd(t_{ab},k)=1$. Otherwise $\gcd(k,t_{ab}')>1$ and hence $\gcd(t_{ab},t_{a})=kd$ in contradiction to the definition of $d$ . Therefore   
\begin{eqnarray*}
1&=&\gcd(t_{ac},b)=\gcd(t_{ac},b+kc),\\
1&=&\gcd(t_{ab},kc)=\gcd(t_{ab},b+kc),\\ 
1&=&\gcd(d,kt_{c}t_{ac})=\gcd(d,dt_{b}t_{ab}'+kt_{c}t_{ac})=\\
 &=& \gcd(ds,dt_{b}t_{ab}' +kt_{c}t_{ac})=\gcd(dsk,dt_{b}t_{ab} +kt_{c}t_{ac})
=\\
 &=&\gcd(t_{a}, dt_{b}t_{ab}'t_{cb}+ kt_{c}t_{ac}t_{cb} )=\gcd(t_{a},b+ kc).
\end{eqnarray*}
\end{proof}
This now allows an unique parametrization of the elements in 
$I_{n}$.

\begin{prop}\label{hilg2}
There is a bijection from the set 
$$P_{n}=\{(c,b):  c\geq 1,\ c\mid n,\quad  b\in \{0,\ldots,\frac{n}{c}-1\},\ 
(c,b)\ n-admissible \}$$ 
to the set $I_{n}$. The map is given by
$$(c,b)\mapsto [c:d(c,b)]$$ with $d(c,b)$ from Definition 
\ref{parametrization}.
\end{prop}

\begin{proof}
We show first that the above map is surjective. For any $[x:y]\in I_{n}$ by 
Proposition \ref{heute} 
there exist an unique $c\geq 1,\ c\mid n$ and $d'$ with $[x:y]=[c:d']$.
 Define $b\in \{0,\ldots,\frac{n}{c}-1\}$ through $d'\equiv 
(b+c)\mod\frac{n}{c} $. 
We claim $(c,b)$ is $n$-admissible. Indeed, from Proposition \ref{gcd} we see  
$\gcd(c,d',n)=1$.
 Assume $\lambda = \gcd(c,b,\frac{n}{c})> 1$. But $\lambda \mid\gcd(d',n)$ and 
hence 
$\lambda \mid \gcd(c,d',n)$. Hence by Lemma \ref{hilg1} there exists $k\in 
\mathbb{Z}$ with 
$\gcd(c,b+k\frac{n}{c})=1.$\\
 Next we claim $d'\equiv d(c,b)\mod\frac{n}{c}$. Indeed 
$d'\equiv(b+c)\mod\frac{n}{c}$ 
and $d(c,b)\equiv(b+c)\mod\frac{n}{c} $ implies $d'=d(c,b)+l\frac{n}{c}$. 
Choose $r,s,t \in \mathbb{Z}$ 
with $rd'-sc-tn=l$. An easy calculation then gives 
$d(c,b)=d'-l\frac{n}{c}=(1-\frac{n}{c}r)d'+(s+t\frac{n}{c})n$ and trivially $c=(1-\frac{n}{c}r)c+rn$. We claim $\gcd((1-\frac{n}{c}r),n)=1$. Obviously $\gcd((1-\frac{n}{c}r),\frac{n}{c})=1$. Assume then  $\gcd((1-\frac{n}{c}r),n)=m>1$. Then $\gcd((1-\frac{n}{c}r),c)=m$. Since 
$d(c,b)=d'-l\frac{n}{c}=(1-\frac{n}{c}r)d'+(s+t\frac{n}{c})n$ the number $m$ 
divides also $d(c,b)$ and 
hence $\gcd(d(c,b),c))>1$ in contradiction to the definition of $d(c,b)$. Hence $[c:d']=[c,d(c,b)].$

To show injectivity of the map $(c,b)\mapsto [c:d(c,b)]$ lets assume $(c',b')$ 
maps to $[c':d(c',b')]$ 
and $[c:d(c,b)]=[c':d(c',b')]$ . Since $c,c' \geq 1$ and $c,c' \mid n$, 
Proposition \ref{heute} shows $c=c'$.
But $b\equiv (d(c,b)-c)\mod\frac{n}{c}$ and $b'\equiv 
(d(c',b')-c)\mod\frac{n}{c}$. 
By Proposition \ref{heute} 
we know that $d(c,b)\equiv d(c',b')\mod\frac{n}{c}$. Therefore also $b\equiv 
b'\mod\frac{n}{c}$ and 
hence $b=b'$ 
since both $b,b' \in \{0,\ldots,\frac{n}{c}-1\}.$
\end{proof}

The set $P_{n}$ can be ordered lexicographically by saying  $(c,b) < (c',b')$ 
iff $c<c'$ or $c=c'$ 
and $b<b'$.

\begin{defi}\label{ABdef}
{\rm
Proposition \ref{hilg2} allows us to identify each element $\pai rs\in I_n$  
with a pair $(c,b)\in P_{n}.$ 
Then we set
\begin{equation}
\label{bekhan}
A_{\pai rs}:=\mat c b 0{\frac{n}{c}},\ 
B_{\pai rs}:=\mat {\frac{n}{c}}0bc
\end{equation}
and define the rational number
$$
x(\pai rs):=\frac{b}{\frac{n}{c}}.
$$ 
Obviously there is a one to one correspondence between the sets $I_{n}$, 
$P_{n}$ and the sets of 
matrices $A\in \Mat_{n}(2,\mathbb{Z})$ which are upper triangular, 
respectively those which are lower 
triangular, and whose entries have greatest common divisor $1$.  
}\hfill$\Box$
\end{defi}


\section{The operator $K$}
\label{opK}

Recall from the introduction the sets of matrices
\begin{eqnarray*}
S_n&=&\left\{\mat abcd : a>c \geq 0, \ d>b\geq 0,\  ad-bc=n\right\}\\
X_n&=&\left\{ \mat c{a}0{\frac{n}{c}} :\ c\mid n,\ 0\leq  a<\frac{n}{c} 
\right\}\\
Y_n&=&\left\{ \mat c0{a}{\frac{n}{c}}: \ c\mid n,\ 0\leq a<c \right\}.
\end{eqnarray*}

\begin{prop}\label{opKiso}
The formula
\begin{equation}\label{Kdef}
K\mat abcd=T^{\sah{\frac{d}{b}}}Q\mat abcd=
 \mat {-c+\sah{\frac{d}{b}}a}{-d+\sah{\frac{d}{b}}b}ab
\end{equation}
defines a bijection $K:S_n\setminus Y_n\to S_n\setminus X_n$ with inverse 
given by the formula
\begin{equation}\label{Kinversedef}
K^{-1}\mat {a'}{b'}{c'}{d'} = MT^{\sah{\frac{a'}{c'}}}QM\mat {a'}{b'}{c'}{d'}=
 \mat {c'}{d'}{-a'+\sah{\frac{a'}{c'}}c'}
{-b'+\sah{\frac{a'}{c'}}d'}.
\end{equation}
\end{prop}

\begin{proof} We denote the right hand side of (\ref{Kdef}) by $\mat 
{a'}{b'}{c'}{d'}$.
The condition $\mat abcd\in S_n\setminus Y_n$ implies
$$a>c\ge 0,\quad d>b>0,\quad ad-bc=n.$$
{}From this it is clear that $c'>0$ so that $\mat {a'}{b'}{c'}{d'}$ is not 
contained in $X_n$.
To show that it is in $S_n$ we  note 
$$a'=\sah{\textstyle\frac{d}{b}}a-c\ge a-c \ge a =c'>0,$$
$$0\le b'=\sah{\textstyle\frac{d}{b}}b-d= 
(\sah{\textstyle\frac{d}{b}}-\textstyle\frac{d}{b})b <b=d',$$
and
$$a'd'-b'c'=(-c+\sah{\textstyle\frac{d}{b}}a)b-(-d+\sah{\textstyle\frac{d}{b}}b)a=ad-bc=n.$$
Thus $K$ is well defined.
That
$$MT^{\sah{\frac{a'}{c'}}}QMT^{\sah{\frac{d}{b}}}Q\mat abcd=\mat abcd$$
follows from 
$$MT^{\sah{r}}QMT^{\sah{s}}Q=\mat 10{r-s}1$$
and $\sah{\frac{a'}{c'}}=\sah{\frac{d}{b}}$ which in turn is a consequence of
$\frac{a'}{c'}=-\frac{c}{a}+\sah{\frac{d}{b}}$ and $-1<-\frac{c}{a}\le 0$.
Similarly we see that
$$MT^{\sah{\frac{d''}{b''}}}QMT^{\sah{\frac{a'}{c'}}}Q\mat 
{a'}{b'}{c'}{d'}=\mat {a'}{b'}{c'}{d'},$$
where $\mat {a''}{b''}{c''}{d''}$ denotes the right hand side of 
(\ref{Kinversedef}).
All that remains to be seen is that $\mat {a''}{b''}{c''}{d''}\in S_n\setminus 
Y_n$ if
$\mat {a'}{b'}{c'}{d'}\in S_n\setminus X_n$, but that can be checked similarly 
as the 
well-definedness of $K$.
\end{proof}
An operator slightly different from the above operator $K$ was used also by 
Choie and Zagier 
in \cite{chza} and by M\"uhlenbruch in \cite{mul} in their derivation of the 
Hecke operators within 
the Eichler, Manin and Shimura theory of period polynomials.
In the following we will use this operator to attach to any index  $ \pai cd 
\in I_{n}$ a sequence of 
elements in $\mathcal{R}_{n}$ which on the other hand are  closely related to 
the minimal partition of 
the rational number $x(\pai cd)$. 

\begin{defi}\label{kidef}{\rm
For $i\in I_n$ we denote by $k_{i}$ the natural number with the property 
that $K^j(A_i)$  (cf. Definition \ref{ABdef}) is  well-defined for $j\leq k_i$ 
and $K^{k_i}(A_i)\in Y_n$. We call
$$A_i,K(A_i),\ldots,K^{k_i}(A_i)$$
the {\it chain} associated with $i\in I_n$.
}\hfill$\Box$ 
\end{defi}

If $A_i\in X_n\cap Y_n$, then clearly $A_i$  forms a chain in itself so that 
$k_i=0$ in this case.

\begin{lemm}
\label{rooz}
Let $\pai cd\in I_n$ and $\part_{x(\pai cd)}=(x_0,x_1,\ldots,x_{k-1},x_k)$
the minimal partition of $x_0=x(\pai cd)=\frac{b}{\frac{n}{c}}$
(cf.~Definition \ref{ABdef} and Remark \ref{minpartrem}).
Suppose that $x_j=\frac{p_j}{q_j}$,  $\gcd(p_j,q_j)=1$, and $q_j\ge 0$.
Then we have $k_{\pai cd}=k-1$ and 
$$K^{j}(A_{\pai cd})=\mat{q_{k-1-j}}{-p_{k-1-j}}{q_{k-j}}{-p_{k-j}}A_{\pai cd}
\quad \forall j=0,\ldots,k-1.$$
\end{lemm}

\begin{proof} Recall the definition of $b\in \{0,\ldots,\frac{n}{c}-1\}$ 
attached to $[c:d]$
in Definition \ref{ABdef}. We assume  $c\geq 1, \quad c\mid n$ and hence 
$$A_{\pai cd}=\mat cb0{\frac{n}{c}}.$$
 We  claim
\begin{equation}\label{Kiter1}
\mat{q_{j-1}}{-p_{j-1}}{q_{j}}{-p_{j}}A_{\pai cd}=
\mat{cq_{j-1}}{bq_{j-1}-\frac{n}{c}p_{j-1}}
{cq_{j}}{bq_j-\frac{n}{c}p_{j}}\in S_n\quad \forall j=1,\ldots,k.
\end{equation}
In fact, using  Lemma~\ref{16jan03}(i) and minimality of the partition we find
$$\left(bq_j-\frac{n}{c}p_{j}\right)-\left(bq_{j-1}-\frac{n}{c}p_{j-1}\right)=
(q_j-q_{j-1})\left(\frac{b}{\frac{n}{c}}-\frac{p_{j-1}-p_j}{q_{j-1}-q_j}\right)>0,$$
whereas
$$\frac{b}{\frac{n}{c}}=x_0\ge x_{j-1} =\frac{p_{j-1}}{q_{j-1}}$$
implies
$bq_{j-1}-\frac{n}{c}p_{j-1}\ge 0$ and even
\begin{equation}\label{notinYn}
bq_{j}-\frac{n}{c}p_{j}> 0\quad \forall j=1,\ldots,k.
\end{equation}
 Since the determinant condition is trivially satisfied
we have proved (\ref{Kiter1}). But there is more detailed information 
available:
Since $\frac{p_0}{q_0}=\frac{b}{\frac{n}{c}}$ and $cq_1< cq_0$  we have 
$$\mat{q_{0}}{-p_{0}}{q_{1}}{-p_{1}}\mat cb0{\frac{n}{c}}=
\mat{cq_{0}}{bq_{0}-\frac{n}{c}p_{0}}
{cq_{1}}{bq_1-\frac{n}{c}p_{1}}=
\mat{cq_{0}}0{cq_{1}}{bq_1-\frac{n}{c}p_{1}}
\in Y_n.
$$
Moreover, 
Remark \ref{gleicheq} shows that 
$\mat {q_{k-1}}{-p_{k-1}}{q_k}{-p_k} =\mat1001$,
so that
$$\mat{q_{k-1}}{-p_{k-1}}0{p_k}\mat cb0{\frac{n}{c}}
=\mat cb0{\frac{n}{c}}\in X_n.$$
On the other hand, by (\ref{notinYn}) none of the 
$\mat{q_{j-1}}{-p_{j-1}}{q_{j}}{-p_{j}}A_{\pai cd}$ in (\ref{Kiter1}) with
$j=1,\ldots,k-1$ can be in $X_n\cup Y_n$ since $q_j\not=0$ for these $j$.
Now it suffices to prove the identities
\begin{equation}
\label{2003miladi}
K:\mat{q_{j}}{-p_{j}}{q_{j+1}}{-p_{j+1}}A_{\pai cd} 
\mapsto \mat{q_{j-1}}{-p_{j-1}}{q_{j}}{-p_{j}}A_{\pai cd}. 
\end{equation}
To prove (\ref{2003miladi}) note first that
for an arbitrary matrix $A\in S_n\setminus Y_n$ we have
$$K(A)A^{-1}=T^{\sah {\frac{d'}{b'}}}Q,$$ 
where $A=\mat {a'}{b'}{c'}{d'}$. 
For $A:=\mat{q_{j}}{-p_{j}}{q_{j+1}}{-p_{j+1}}A_{\pai cd}$ by 
Lemma~\ref{16jan03}(ii) 
we have
$$\sah {\textstyle\frac{d'}{b'}}
=\sah {\frac{bq_{j+1}-\frac{n}{c}p_{j+1}}{bq_j-\frac{n}{c}p_{j}}}
=\sah {\frac{xq_{j+1}-p_{j+1}}{xq_j-p_{j}}}
=p_{j-1}q_{j+1}-p_{j+1}q_{j-1}
$$
and calculate
\begin{eqnarray*}
&&\hskip -4em \mat{q_{j-1}}{-p_{j-1}}{q_{j}}{-p_{j}}A_{\pai cd}
(\mat{q_{j}}{-p_{j}}{q_{j+1}}{-p_{j+1}}A_{\pai cd})^{-1}=\\
&=&\mat{q_{j-1}}{-p_{j-1}}{q_{j}}{-p_{j}}\mat{-p_{j+1}}{p_{j}}{-q_{j+1}}{q_{j}}\\
&=&\mat {p_{j-1}q_{j+1}-p_{j+1}q_{j-1}}{-1}{1}{0}\\
&=&T^{(p_{j-1}q_{j+1}-p_{j+1}q_{j-1})}Q\\
&=&T^{\sah {\frac{d'}{b'}}}Q
\end{eqnarray*}
proving
$$K(A)=\mat{q_{j-1}}{-p_{j-1}}{q_{j}}{-p_{j}}A_{\pai cd}.$$  
\end{proof}
\begin{rem}
{\rm
A construction rather similar to the one in Lemma \ref{rooz} has been used also by L. Merel in \cite{mer}, where he  discussed the connection between the ordinary Hecke operators for the group $\Gamma_{0}(n)$ and continued fractions.
\hfill$\Box$}
\end{rem}


\section{The Lewis equations for the group $\Gamma_{0}(n)$}
\label{lew}

Consider the right $\R$-module 
$$\R^{I_n}:=\{\psi\colon I_n\to \R\}=\R\otimes \Z^{I_n},$$
whose elements we also denote by $\psi=(\psi_i)_{i\in I_n}$.
The module $\R^{I_n}$ is equipped with a natural left $\GL(2,\Z)$-action
given by $(g\cdot \psi)_i=\psi_{ig}$. Recall  the right regular representation
$\rho$ of $\GL(2,\Z)$ on $\ga n\backslash \GL(2,\Z)$ which we identify with 
$I_n$.
Then we have
$$g\cdot(R\otimes w)=R\otimes \rho(g)w\quad \forall g\in \GL(2,\Z), w\in 
\Z^{I_n}, R\in \R$$
and, by abuse of notation we write $\rho(g)\psi$ for $g\cdot \psi$.
It is important to note that
$$(\rho(g)\psi)R=\rho(g)(\psi R)\quad\forall g\in \GL(2,\Z), \psi\in \R^{I_n}, 
R\in \R.$$
Therefore all terms in the equation 
\begin{equation}
\label{puf}
\psi-\rho(T^{-1})\psi T-\lambda \rho(T^{-1}M)\psi TM \equiv 0 \mod 
(\I^\lambda)^{I_n}
\end{equation} 
for the unknown $\psi\in \R^{I_n}$  are unambiguous.
Comparing (\ref{puf}) with equation (\ref{27JAN}) shows that it is reasonable 
to 
call
(\ref{puf}) the {\it Lewis equation} in $(\I^\lambda\backslash\R)^{I_n}$ 
corresponding to the transfer operator for the group $\Gamma_{0}(n)$.

\begin{rem}\label{scalLewisinRIn}{\rm 
To rewrite  (\ref{puf}) as a system of scalar equations
we have to determine the actions of $T^{-1}$ and $T^{-1}M$ on
the indexing coset space $I_n$. They are given by 
$$
\pai cd T^{-1} = \pai c{d-c}
\quad \mbox{and} \quad
\pai cd T^{-1}M= \pai {d-c}c.
$$
Thus the corresponding system of scalar equations is
\begin{equation}
\label{18nov02}
\psi_{\pai cd}-\psi_{\pai c{d-c}}T-\lambda\psi_{\pai {d-c}c}TM\equiv 0 \mod 
\I^\lambda,\quad
\forall \pai cd\in I_n.  
\end{equation}
Replacing $\pai cd$ by $\pai cd\mat {-1}011=\pai {d-c}d$, 
 multiplying the resulting equation from the right by $\lambda M$, and then 
subtracting it 
from the original equation we get 
\begin{equation}
\label{zagier}
\psi_{\pai cd}\equiv \lambda\psi_{\pai {d-c}d}M \mod \I^\lambda
\end{equation}
We call $\psi_{\pai cd}, \psi_{\pai {d-c}d}$ a {\it symmetric pair}. 
It is easy to see that ${\pai cd}={\pai {d-c}d}$ if and only if
$c=1$ and $n\mid d(d-2)$.
In particular,  $\pai 10$ and $\pai 12$ have this property. We call an
element $\psi_{\pai cd}$ with this property  a {\it self-symmetric 
element}.
}\hfill$\Box$
\end{rem}

There is a close relation between the equations 
(\ref{18nov02}) and (\ref{27JAN}). To explain this relation
we have to write (\ref{27JAN}) as a 
system of scalar equations.

\begin{rem} \label{scalLewisGaussKuzmin}{\rm
Let $\Phi\in \mathcal{B}(D)\otimes \C^{I_n}$ be a solution of (\ref{27JAN}).
We write it as $\Phi=(\phi_i)_{i\in I_n}$ with $\phi_i\in \mathcal{B}(D)$. 
Similarly  
as in Remark \ref{scalLewisinRIn} we see that
the $\phi_i$ satisfy the following system of scalar equations 
\begin{equation}
\label{7feb03}
\phi_{[c:d]}-\phi_{[c:d-c]} T-\lambda \phi_{[d-c:c]} TM =0 \quad \forall 
[c:d]\in I_n.
\end{equation}
}\hfill$\Box$
\end{rem}

\begin{rem}\label{LewisRelation}{\rm
Let $\sol$ be a solution of the scalar 
Lewis equation (\ref{2feb03}). Then according to \cite{leza2} we have $\phi\in 
\F_0$
so that in view of Remark \ref{G+F0} equation (\ref{2feb03}) can be rewritten 
as
\begin{equation}\label{scallewis2}
\sol\mid_s{(I-T-\lambda TM)}=0
\end{equation}
Next we write equation (\ref{18nov02}) as 
\begin{equation}
\label{26nov02}
\psi_{[c:d]}-\psi_{[c:d-c]} T-\lambda \psi_{[d-c:c]} TM =(I-T-\lambda 
TM)P_{[c:d]},
\end{equation}
with  $P_{[c:d]}$ some element in $\mathcal{R}$.
Now we want both sides of (\ref{26nov02}) to act on $\phi$ via the 
slash action. If this were possible, by (\ref{scallewis2}) 
the right hand side would annihilate $\phi$, so that
$\Phi:=(\phi\mid_s \psi_i)_{i\in I_n}$ in view of
Remark \ref{scalLewisGaussKuzmin} were a solution of (\ref{27JAN}).
Lemma \ref{2bahman81} shows that this is indeed possible as long as all the 
matrices occurring
in $P_{[c:d]}$ satisfy the hypotheses of this lemma.
}
\hfill$\Box$
\end{rem}

\begin{rem}\label{IlambdaI}{\rm
For $\lambda = \pm 1$ we have 
\begin{equation}\label{lambdaI-M}
\lambda I-M=(I-T-\lambda TM)(\lambda I-M) \equiv 0\mod \I^\lambda
\end{equation}
and so
$$
(I-T-MTM)=(I-T-\lambda TM) + (\lambda I-M)TM=(I-T-\lambda TM)(I+\lambda 
TM-MTM)
$$
which implies that $\I:=(I-T-MTM)\R\subset\I^\lambda$.
Replacing $\pai cd$ by $\pai c{d+c}=\pai cd \mat 1101$ in (\ref{18nov02}) 
we arrive at the system
\begin{equation}
\label{2003}
\psi_{\pai c{d+c}}-\psi_{\pai cd}T-\lambda\psi_{\pai dc}TM \equiv 0 \mod 
\I^\lambda \quad \forall 
\pai cd\in I_n
\end{equation}
of scalar equations. A further modification is suggested by (\ref{lambdaI-M}):
\begin{equation}
\label{2003b}
\psi_{\pai c{d+c}}-\psi_{\pai cd}T-M\psi_{\pai dc}TM \equiv 0 \mod \I\quad 
\forall 
\pai cd\in I_n.
\end{equation}
This is the equation we will solve and from where we will construct a 
solution of (\ref{2003}). More precisely, if
$\psi_{\pai c{d+c}}-\psi_{\pai cd}T-M\psi_{\pai dc}TM=(I-T-MTM)R,\ R\in\R$,
then
\begin{eqnarray*}
&&\hskip -2em \psi_{\pai c{d+c}}-\psi_{\pai cd}T-\lambda\psi_{\pai dc}TM=\\
&=&(I-T-\lambda TM)((1+\lambda TM-MTM)R+(M-\lambda I)\psi_{\pai dc}TM)\\
&\equiv& 0 \ \mod \ \I^\lambda .
\end{eqnarray*}
In view of Remark \ref{gula} 
if  $R\in\Z[\T]$ and $\psi_{\pai dc}\in \Z[\Mat_*(2,\Z^+\cup\{0\})]$ then
$$
(1+\lambda TM-MTM)R+(M-\lambda I)\psi_{\pai dc}TM\in \Z[\T]
$$
 This means that if for $\pai cd\in I_n$ we have that  $\psi_{\pai cd}\in 
\Z[\Mat_*(2,\Z^+\cup\{0\})]$ 
are solutions 
of (\ref{2003b}) fulfilling the condition of  Lemma \ref{2bahman81},
then these $\psi_{\pai cd},\ \pai cd\in I_n$ solve  also (\ref{2003}) and 
satisfy
the condition of  Lemma \ref{2bahman81}.
}\hfill$\Box$
\end{rem}

\begin{lemm}
\label{mayer}
For any two partitions $\part_1,\part_2$ of $x\in \Q^+$ we have
$$m(\part_1)- m(\part_2) \equiv 0 \mod \I.$$ 
If  $s\in \C$ and $\sol$ is a solution of the scalar Lewis equation
(\ref{2feb03}), then we have the following (well-defined) equality
$$
0=\sol\mid_s \big(m(\part_1)-m(\part_2)\big).
$$
\end{lemm}

\begin{proof}
By Lemma ~\ref{ezdev} it is enough to prove the lemma for a 
partition $\part$ and its modification
$\part(l)$ in (\ref{13jan03}). In the notation of (\ref{13jan03}) the element
 $m\big(\part(l)\big)$ is given by
$$
\ldots+\mat{q_{l-1}}{-p_{l-1}}{q_l+q_{l-1}}{-p_l-p_{l-1}}+
\mat{q_{l}+q_{l-1}}{-p_{l}-p_{l-1}}{q_{l}}{-p_{l}}+
\mat{q_l}{-p_l}{q_{l+1}}{-p_{l+1}}\ldots
$$
so that
\begin{eqnarray*}
&&\hskip -2em m\big(\part(l)\big)-m(\part)=\\
&=&
\mat{q_{l-1}}{-p_{l-1}}{q_l+q_{l-1}}{-p_l-p_{l-1}}+
\mat{q_{l}+q_{l-1}}{-p_{l}-p_{l-1}}{q_{l}}{-p_{l}}-
\mat{q_{l-1}}{-p_{l-1}}{q_l}{-p_l}\\
&=&(MTM+T-I)\mat{q_{l-1}}{-p_{l-1}}{q_l}{-p_l}.
\end{eqnarray*}
If $q_l=0$, then we have $-q_{l-1}p_l=1$ so that $q_{l-1}>0$ and $-p_l>0$,
which in turn shows that $\mat{q_{l-1}}{-p_{l-1}}{q_l}{-p_l}$ satisfies the 
conditions of
Lemma \ref{2bahman81}. 
If $q_l>0$ and $q_{l-1}=0$, then $q_lp_{l-1}=1$ so that $p_{l-1}>0$
and $x_l=\infty$. This contradiction shows that $q_{l-1}>0$ and hence also in 
this case
 $\mat{q_{l-1}}{-p_{l-1}}{q_l}{-p_l}$ satisfies the conditions of
Lemma \ref{2bahman81}. 
As a result we obtain that
$\phi\mid_s \big(m\big(\part(l)\big)-m(\part)\big)$ is well-defined and in 
fact equal to $0$ since
$\phi$ is a solution of (\ref{2feb03})
\end{proof}

\begin{prop}\label{psim(P)}
For $[c:d]\in I_n$ set 
$$\widetilde\psi_{[c:d]}:=\sum_{j=0}^{k_{[c:d]}}K^j(A_{[c:d]}).$$
\begin{enumerate}
\item[{\rm(i)}] 
$\widetilde\psi\in \R^{I_n}$ is a solution of (\ref{2003b}).
\item[{\rm(ii)}] If  $s\in \C$ and $\sol$ is a solution of the scalar Lewis 
equation
(\ref{2feb03}), then 
$$0=\phi\mid_s\big(\widetilde \psi_{\pai cd}T+M\widetilde\psi_{\pai dc}TM
-\widetilde\psi_{\pai {c}{c+d}}\big).$$ 
\end{enumerate}
\end{prop}

\begin{proof}
Note first that Lemma \ref{rooz}, together with the   Definitions \ref{ABdef} 
and  \ref{mdef}, shows
\begin{equation}\label{tildepsiformel}
\widetilde \psi_{[c:d]}=m(P_{x([c:d])})A_{[c:d]}.
\end{equation}

Fix $[c:d]\in I_n$. According to Lemma ~\ref{heute} we may assume that 
$d=ic'$ with $1\le c,c'\mid n$ and $\gcd(c,c')=1=\gcd(i,n)$.
Thus we are in the situation of Lemma \ref{XTslemma}.
But in the notation of that lemma we have
\begin{eqnarray*}
&&\hskip -4em\widetilde \psi_{\pai cd}T+M\widetilde\psi_{\pai dc}TM=\\
&=&m(\part_{x}) \mat c{(ic'-c)_{\frac{n}{c}}}0{\frac{n}{c}}T
+
M\, m(\part_{z})
\mat {c'}{(\hat \imath c-c')_{\frac{n}{c'}}}0{\frac{n}{c'}}TM\\
&=&m\big((\part_{z}\cdot X)\vee (\part_x\cdot T^{s})\big)
\mat c{(ic')_{\frac{n}{c}}}0{\frac{n}{c}}.
\end{eqnarray*}
On the other hand we have
$$
\widetilde\psi_{\pai {c}{c+d}}
=m(\part_{y})\mat c{(d)_{\frac{n}{c}}}0{\frac{n}{c}}
$$ 
so that now the first part of Lemma \ref{mayer} shows that
$$\widetilde \psi_{\pai cd}T+M\widetilde\psi_{\pai dc}TM
-\widetilde\psi_{\pai {c}{c+d}}\equiv 0\mod \I.$$
This proves (i) and (ii) is now an immediate consequence of the second part of 
Lemma \ref{mayer}. 
\end{proof}

\begin{rem}{\rm 
Replacing $\pai cd$ with $\pai c{d-c}=\pai cd $ 
in (\ref{2003b}) we find 
\begin{equation}
\label{2003c}
\psi_{\pai cd}-\psi_{\pai c{d-c}}T-M\psi_{\pai {d-c}c}TM \equiv 0 \mod \I\quad 
\forall 
\pai cd\in I_n
\end{equation}
and by Remark \ref{IlambdaI} a solution of (\ref{2003c}) is also a solution of
(\ref{26nov02}). Thus Proposition \ref{psim(P)} actually provides a solution
of (\ref{26nov02}). Moreover, together with the last part of  Remark 
\ref{IlambdaI}
and Remark \ref{LewisRelation}, it provides the equality 
$$\phi\mid_s\left(\widetilde \psi_{\pai cd}T+\lambda\widetilde\psi_{\pai dc}TM
-\widetilde\psi_{\pai {c}{c+d}}\right) = 0$$
so that 
$\Phi:=(\phi\mid_s \psi_i)_{i\in I_n}$ is a solution of (\ref{27JAN}).
Thus we have now proved Theorem \ref{haupt}.
\hfill$\Box$}
\end{rem}


\section{Hecke operators}
\label{Heckeop}

In \cite{zag1}, \cite{zag2} D. Zagier derived a representation of the Hecke operators on the 
space of period polynomials 
for the group $\PSL (2,\mathbb{Z})$ by transferring the action of the 
classical Hecke operators on the 
space of cusp forms via the Eichler-Shimura-Manin isomorphism to the space of 
period polynomials.
 In his thesis T. M\"uhlenbruch (\cite{mul}) found another representation for 
these operators in terms of 
matrices with nonnegative entries 
which allowed him to extend their action to the space of period functions with 
arbitrary weight. 
It turns out that the special solutions of the Lewis equations for the 
congruence subgroups $\Gamma_{0}(n)$ 
we constructed in Proposition \ref{psim(P)} are closely related to the Hecke 
operators for 
$\PSL (2,\mathbb{Z})$ in the form given by  M\"uhlenbruch.

Indeed, since both the maps $T:I_{n}\to I_{n}$ and $MT:I_{n}\to I_{n}$ are 
invertible any solution 
$\Phi$ of the Lewis equation (\ref{DFG}) for $\Gamma_{0}(n)$ given by  
$\sol_i=\sol_i(z),i\in I_{n}$,  
determines a solution $\tilde{\phi}=\tilde{\phi}(z)$ of the Lewis equation 
(\ref{2feb03}) for the group 
$\PSL (2,\mathbb{Z})$ with
\begin{equation}\label{Sum}
\tilde{\phi}(z)=\sum_{i\in I_{n}}\sol_i(z).
\end{equation}
Clearly, it can happen that this function vanishes identically. This just 
signals that the corresponding 
solution $\sol_i, i\in I_{n}$ for the group  $\Gamma_{0}(n)$ is not related to 
any solution $\phi$ of the 
group $\PSL (2,\mathbb{Z})$ and hence, in analogy to the Atkin-Lehner theory, 
should be called  a 
{\it new solution} of the Lewis equation for $\Gamma_{0}(n)$.
 The special solution, however, determined in Proposition \ref{psim(P)} leads 
to a nontrivial solution 
$\tilde{\phi}$ which furthermore depends linearly on the solution $\phi$ of 
equation (\ref{2feb03}).
 This shows that the map $\tilde{H}_{n}: \phi \mapsto \tilde{\phi}$ with 
$\tilde{\phi}$ as defined in 
equation (\ref{Sum})
 determines a linear operator in the space of period functions of the group 
$\PSL (2,\mathbb{Z})$.
To determine the explicit form of the operator $\tilde{H}_{n}$
 we have to characterize the matrices $K^{j}(A_{[c:d]})$ appearing in the 
definition of the solutions 
$\tilde{\psi}_{[c:d]}$ in Proposition \ref{psim(P)} in more detail.

From the definition of the operator $K$ in Proposition \ref{opKiso} it is 
obvious that all matrix elements 
of $A \in S_{n}\setminus Y_{n}$ have greatest common divisor $1$ 
if and only if the matrix elements of the matrix  $K(A)$ have this property. 
Since the entries of the matrix $A_{[c:d]}$ in Definition \ref{ABdef} for 
$[c:d]\in I_{n}$  have greatest 
common divisor $1$ all the matrices appearing in the definition of 
$\tilde{\psi}_{[c:d]}$ in 
Proposition \ref{psim(P)} have this property. 

Consider next any matrix $A\in S_{n}\setminus X_{n}$ whose entries have 
greatest common divisor $1$. 
If $A=\mat abcd$, then $K^{-1}A=\mat {a'}{b'}{c'}{d'}$ with $c'<c$ and hence 
there exists $j\in \mathbb{N}$ 
with $K^{-j}A \in X_{n}$. 
But from Proposition \ref{hilg2} it follows that any matrix $A$ in $X_{n}$ 
whose entries have only $1$ as 
a common divisor  appears as  $A_{[c:d]}$ for some $[c:d]\in I_{n}$. 
This shows that any matrix $A$ in the set $S_{n}$ whose entries have no common 
divisor besides $1$ appears 
exactly once in one of the components $\tilde{\psi}_{[c:d]}$ in Proposition 
\ref{psim(P)}.

Denote then by $\tilde{T}_{n}$ the matrix
$$\tilde{T}_{n}:= \sum_{A \in S_{n}:\gcd(a,b,c,d)=1}A, \quad A=\mat abcd .$$ 
Then one finds for the operator $\tilde{H}_{n}$ acting on the space of period 
functions $\phi$ for the 
group $\PSL (2,\mathbb{Z})$ $$\tilde{H}_{n} \phi=\phi\mid_{s}\tilde{T}_{n}$$
Summarizing we have shown:

\begin{theo}\label{Hecke} For any solution $\phi = \phi(z)$ of the Lewis 
equation (\ref{2feb03}) for 
$\PSL (2,\mathbb{Z})$ with arbitrary weight $s$ 
the function $\tilde{\phi}=\tilde{\phi}(z)= \tilde{H}_{n}\phi 
(z)=\phi\mid_{s}\tilde{T}_{n}(z)$ is also a solution of 
equation (\ref{2feb03}) with weight $s$.
\end{theo}

Comparing the operators $\tilde{T}_{n}$ with the Hecke operators $T_{n}$ of 
M\"uhlenbruch and Zagier in 
(\ref{madar}) we find as a corollary
\begin{cor}
The operators $\tilde{T}_{n}$ and the Hecke operators $T_{n}$ defined in (2) 
are related through
$$T_{n}=\sum_{d^{2}\mid n}\mat d00d \tilde{T}_{\frac{n}{d^{2}}}.$$ The 
operators coincide if and only if 
$n$ is a product of distinct primes.
\end{cor}

The operators $\tilde{T}_{n}$ have been constructed from special solutions of 
the Lewis equation (3) 
for the group $\Gamma_{0}(n)$. It turns out that also the Hecke operators 
$T_{n}$ can be derived in 
this way. To do this consider any $n_1,n_2\in\mathbb{N}$.  
For $n_2\mid n_1$ there is a canonical surjective map 
$\sigma_{n_1,n_2}:[\Z\times\Z]_{n_1}\rightarrow [\Z\times\Z]_{n_2}$ which 
is equivariant with respect to the $\GL(2,\Z)$-actions, i.e.
\begin{equation}
\label{vietnam}
\sigma_{n_1,n_2} ([x:y])A=\sigma_{n_1,n_2} ([x:y] A),
\quad \forall A\in \GL(2,\Z),\ [x:y]\in [\Z\times\Z]_{n_1}
\end{equation}
Therefore (cf. Remark \ref{Ininjection})  $\sigma_{n_1,n_2}$
induces a map $I_{n_1}\rightarrow I_{n_2}$ which we still denote by 
$\sigma_{n_1,n_2}$ (or simply by $\sigma$ if $n_1$ and $n_2$ are clear from 
the 
context). 

\begin{prop}
\label{bekhanid}
If $\psi_i^1,\ i\in I_{n_1}$ solve  ~(\ref{18nov02}) for $n=n_1$ then
$\psi_j^2:=\sum_{i\in\sigma^{-1}j} \psi_i^1,\  j\in I_{n_2}$ solve  
~(\ref{18nov02}) for $n=n_2$. Moreover, if $\psi_j^2,\ j\in I_{n_2}$ solve 
 ~(\ref{18nov02}) for  $n=n_2$ then $\psi_i^1:=\psi_{\sigma i}^2,\ i\in 
I_{n_1}$ 
solve  ~(\ref{18nov02}) for $n=n_1$.  
 \end{prop} 

\begin{proof}
~(\ref{vietnam}) implies that the fibers $\sigma^{-1}j$ of $\sigma$ 
are invariant under the action of $\GL(2,\Z)$, and in particular
$T^{-1}$ and  $T^{-1}M$. 
\end{proof}

Proposition \ref{bekhanid} shows that any solution $\Phi= (\phi_{i}, i\in 
I_{n})$ of equation 
(\ref{18nov02}) for the group $\Gamma_{0}(\frac{n}{d^{2}})$ determines a 
solution for this equation for 
the group $\Gamma_{0}(n)$ whose components coincide with the components for 
the former group. 
Indeed any component shows up $\mu$-times, where $\mu$
is the index of $\Gamma_{0}(\frac{n}{d^{2}})$ in $\Gamma_{0}(n)$.   
 Taking for $\Phi$ the special solution $\phi\mid_{s}\tilde{\psi}_{i},i\in 
I_{\frac{n}{d^{2}}}$ 
determined in (ii) of Proposition \ref{psim(P)} we therefore get

\begin{cor}\label{Teiler}
For any solution $\phi$ of the Lewis equation (\ref{2feb03}) for the group 
$\PSL (2,\mathbb{Z})$ with 
weight $s$ the functions 
$\tilde{\phi}_{j,d}:=\phi\mid_{s}\tilde{\psi}_{\sigma(j)},j\in I_{n}$ define 
a 
solution of the Lewis equation (3) for the group $\Gamma_{0}(n)$ with weight 
$s$.
\end{cor}

Hence also the function $\tilde{\phi}_{d}$ with
$$\tilde{\phi_{d}}=\frac{1}{\mu}\sum_{j\in I_{n}}\tilde{\phi}_{j,d}
= \phi\mid_{s}\tilde{T}_{\frac{n}{d^{2}}}$$ 
defines a solution of the Lewis equation (\ref{2feb03}) for the group $\PSL 
(2,\mathbb{Z})$. 
Obviously the matrix inducing this solution $\tilde{\phi}_{d}$ coincides with 
the matrix 
$\tilde{T}_{\frac{n}{d^{2}}}$. This shows that indeed the Hecke operator $T_n$ 
on the period functions of 
$\PSL (2,\mathbb{Z})$ for arbitrary weight $s$ can be derived from special 
solutions of the Lewis equation 
for the group $\Gamma_{0}(n)$ with weight $s$.

The extension of this approach to the Hecke operators on period functions for 
the congruence subgroups 
$\Gamma_{0}(n)$ will be discussed in a forthcoming paper.\\
\\

{\bf Acknowledgements}

This work has been supported by the Deutsche Forschungsgemeinschaft through the DFG Forschergruppe ``Zetafunktionen und lokalsymmetrische R\"aume''. One of the authors (D.M.) thanks the IHES for financial support and the kind hospitality extended to him during the final preparation of this paper.

\end{document}